\documentclass[noamsfonts,reqno,11pt]{amsart}
\usepackage[foot]{amsaddr}         

\usepackage[sfdefault,mono,vvarbb]{notomath} 

\usepackage[svgnames]{xcolor}
\usepackage{graphicx}

\usepackage[labelfont=bf]{caption}
\usepackage{subcaption}
\usepackage{nicematrix}
\usepackage{booktabs}

\usepackage[shortlabels]{enumitem}
\setlist[description]{leftmargin=\parindent,labelindent=\parindent}
\usepackage{hyperref}

\usepackage{bm}
\usepackage{accents}

\usepackage{algorithm}
\usepackage{algpseudocode}

\algblock[WHERE]{Where}{EndWhere}
\algblockdefx[WHERE]{Where}{EndWhere}%
   [1]{\textbf{where} #1}%
   {\textbf{end where}}
\algcblockdefx[ELSEWHERE]{WHERE}{ElseWhere}{EndWhere}%
   [1][]{\textbf{else where} #1}
   {\textbf{end where}}

\usepackage[numbers,sort&compress]{natbib} 

\calclayout

\input{def}

\begin{document}

\title%
[
Adaptive multigrid for high-order DG methods
]
{
Adaptive multigrid for high-order discontinuous Galerkin methods
based on the full approximation scheme
}
\author{Jörg Stiller}
\address{%
  TU Dresden, Institute of Fluid Mechanics, 01062 Dresden, Germany}
\email{joerg.stiller@tu-dresden.de}

\begin{abstract}
We propose an adaptive multigrid (MG) method for discontinuous Galerkin formulations of elliptic problems using Brandt’s full approximation scheme (FAS). Unlike common approaches, this method achieves local $\mathit{hp}$-refinement of hexahedral meshes without the need for hanging nodes. The core component of the FAS-MG method is an overlapping Schwarz smoother, which is optionally accelerated by a Krylov method. This smoother is designed for unstructured curvilinear meshes but maintains a tensor-product structure for fast diagonalization. Numerical experiments demonstrate the exceptional efficiency of the FAS-MG method. Dedicated studies confirm its robustness against high aspect ratios, element deformation, and irregular mesh topology. We also verify its capability for dynamic parallel mesh adaptation using the wave-front benchmark of \v{C}erven\'{y}, Dobrev, and Kolev (SIAM J. Sci. Comp. 41, 2019). Finally, we present preliminary results of extending the method to incompressible Navier-Stokes problems.
\end{abstract}


\maketitle



\section{Introduction}
\label{sec:introduction}

High-order discretization methods, adaptive mesh refinement, and multigrid
solution methods are essential components in the development of efficient
numerical models capable of addressing future challenges in fluid dynamics,
meteorology, astrophysics, and other research areas based on partial
differential equations.

In computational physics, higher-order discretization methods are widely used 
and are increasingly gaining importance.
The most important approaches include 
spectral methods \cite{SE_Canuto2011a}, 
finite element methods
\cite{FE_Ern2004a,SE_John2016a}, 
spectral element methods \cite{SE_Deville2002a,SE_Hesthaven2008a,SE_Karniadakis2005a}, 
as well as 
finite difference and finite volume methods \cite{SE_Shu2003a}. 
Here, we use the discontinuous Galerkin spectral-element method (DG-SEM).
However, the focus is on adaptive multigrid methods for high-order methods,
which are briefly reviewed below.

\subsection{Multigrid for high-order FEM and SEM}

The mathematical foundations of multigrid methods were established in the
1970s and 1980s. 
Early work focused on low-order finite-difference, finite-volume and finite-element 
methods \cite{MG_Brandt1977a,MG_Hackbusch1985a,MG_McCormick1986a}.
To be efficient with high-order finite element methods (FEM) or spectral
element methods (SEM), multigrid methods must account for the strong coupling
of degrees of freedom inside the elements and exploit the structure of the
discrete operators. A first milestone was achieved by \citet{MG_Lottes2005a},
who developed polynomial multigrid methods with overlapping Schwarz smoothers
for continuous Galerkin spectral-element methods (CG-SEM). On isotropic
Cartesian grids, these methods reduced the residual by ten orders of magnitude
within approximately ten iterations. Their unprecedented computational
efficiency resulted from exploiting sum factorization and fast diagonalization \cite{MG_Lynch1964a} for the evaluation and inversion of tensor-product operators. Acceleration with conjugate-gradient and GMRES
techniques improved robustness against high aspect ratios. The method was
subsequently extended to unstructured curvilinear meshes and applied to
incompressible Navier-Stokes problems \cite{MG_Fischer2005a}.

In \cite{MG_Stiller2016a}, the method was improved by introducing a polynomial
weighting function for combining the corrections from overlapping subdomains.
This modification reduces the iteration count by a factor of up to 3 and
yields runtime savings of approximately 50 percent. Further substantial
improvements can be achieved using static condensation
\cite{MG_Haupt2013a,SE_Huismann2019a}. In \cite{MG_Stiller2016b}, the
Schwarz multigrid method was generalized to DG-SEM on two- and
three-dimensional Cartesian grids.

For unstructured meshes, non-overlapping block preconditioners have been
investigated, for example, in \cite{MG_Kanschat2004a,MG_Bastian2019a}.
However, these approaches become less efficient as the polynomial increases
and are not competitive with overlapping Schwarz multigrid methods
\cite{MG_Stiller2016b}. \citet{MG_Witte2020a} proposed overlapping Schwarz 
smoothers based on vertex-centered subdomains comprising whole elements. 
The large overlap yields exceptional multigrid convergence rates but
also incurs high computational costs, particularly on unstructured meshes,
where fast diagonalization is not feasible.
\citet{MG_Vincent2019a} and \citet{MG_Vu2022a} developed a Schwarz multigrid
method for DG-SEM based on overlapping element-centered subdomains including 
only face neighbors.
This approach simplifies the handling of unstructured meshes.
However, it also prevents the use of fast diagonalization and disregards
potentially valuable information.
In addition to Schwarz methods, Chebyshev methods have been used as 
smoothers in multigrid methods for high-order FEM and SEM \cite{MG_Sundar2015a,MG_Kronbichler2018a}. More
recently, Chebyshev and other polynomial methods have been employed to
accelerate Schwarz multigrid methods \cite{MG_Phillips2025a}.

\subsection{Adaptive multigrid}

Adaptive techniques require multigrid methods suitable for locally refined
meshes. This issue was already addressed in early work, for example, in
\cite{MG_Brandt1977a,MG_Hackbusch1985a,MG_McCormick1986a}. Most methods seek
the solution on a locally refined grid covering the entire computational
domain. The resulting multilevel problem is typically solved using 
either global coarsening
\cite{MG_Becker2000a,MG_Mitchell2010a,MG_Vincent2019a,MG_Vu2022a,
MG_Munch2023a,MG_Mann2025a} or local smoothing
\cite{MG_Kanschat2004a,MG_Janssen2011a,MG_Fortunato2019a,MG_Clevenger2020a}.
In the former case, the smoother sweeps over the entire domain on each level,
whereas local smoothing restricts the sweeps to regions refined on the
corresponding level.

Using a global mesh on the finest level provides the adapted solution in a
single piece and simplifies the implementation of conservative methods. On
the other hand, it requires transition elements or hanging nodes. Transition
elements preclude the use of purely hexahedral meshes, whereas hanging nodes
require local modifications of the discrete operators. These complications
can be avoided by adopting the full-approximation storage (FAS) algorithm of
\citet{MG_Brandt1977a}. Its central idea is to confine the meshes to the
respective refinement zones and store the full solution instead of the
correction on each level. Consecutive levels are coupled by transfer operators
without requiring hanging nodes. FAS integrates seamlessly with the full
multigrid (FMG) method, which provides an efficient startup procedure
\cite{MG_Brandt1977a}.

Adaptive FAS multigrid methods have been developed, for example, in
\cite{MG_Bai1987a,MG_Thompson1989a,MG_Lee2007a,MG_Feng2018a}. All of these
studies used second-order finite-difference or finite-volume methods for
discretization. Apart from his own work on space-time adaptive
DG-SEM for one-dimensional conservation laws \cite{TI_Pfister2026b}, the
author is not aware of any applications of FAS multigrid to high-order
adaptive methods.

\subsection{Contribution of the present paper}

This paper develops an efficient and robust multigrid method for high-order
discontinuous Galerkin methods on locally refined meshes. The method adapts
the full-approximation storage algorithm to DG-SEM. Its main ingredient is a 
weighted overlapping Schwarz smoother that extends prior work
\cite{MG_Stiller2016b,MG_Stiller2017a} to unstructured curvilinear hexahedral
meshes. In contrast to the approach proposed in
\cite{MG_Vincent2019a,MG_Vu2022a}, the smoother also includes data from edge
and vertex neighbors while maintaining the tensor-product structure required
for fast diagonalization. Acceleration using the inexact conjugate-gradient
method \cite{KR_Golub1999a} renders the smoother exceptionally robust against
high aspect ratios, element deformation, and irregular mesh topology. In
addition, hierarchical start strategies such as full multigrid are adopted to
further accelerate the solution process. Numerical experiments demonstrate
excellent efficiency with global or local $\mathit{hp}$-refinement and 
competitiveness with other approaches.

The remainder of the paper is organized as follows: Section~\ref{sec:dg-method}
describes the elliptic model problem and its discretization with DG-SEM.
Section~\ref{sec:ml-formulation} introduces the multilevel mesh and derives the
multilevel DG-SEM formulation based on the FAS scheme. Building on this,
Section~\ref{sec:mg-method} develops the multigrid solver.
Section~\ref{sec:implementation} gives a short overview of the data structures,
mesh adaptation and partitioning, and the implementation of the method.
Section~\ref{sec:num-exp} presents numerical experiments that assess
efficiency and robustness and demonstrate the capability for local and
adaptive mesh refinement. Section~\ref{sec:outlook} gives an outlook on the
application to flow problems. Finally, Section~\ref{sec:conclusion} summarizes
the findings and identifies perspectives for future research.



\section{Discontinuous Galerkin method}
\label{sec:dg-method}


\subsection{Problem definition}
\label{sec:problem}

As a model problem we consider the elliptic PDE
\begin{subequations}
\label{eq:helmholtz}
\begin{equation}
  \label{eq:helmholtz:pde}
  \kappa u - \nabla \cdot \nu \nabla u = f
  \quad \text{in} \quad \Omega
\end{equation}
with boundary conditions
\begin{align}
  \label{eq:helmholtz:bc:dirichlet}
  u &= u_{\mathsc d} \quad \text{on} \quad \d\Omega^{\mathsc d} \\
  \label{eq:helmholtz:bc:neumann}
  \nu \d_n u &= q_{\mathsc n} \quad \text{on} \quad \d\Omega^{\mathsc n}
\end{align}
\end{subequations}
where
${\Omega \subset \mathbb R^3}$
is a possibly periodic, bounded domain,
${\kappa \ge 0}$ 
is a non-negative constant,
${\nu(\V x) > 0}$ 
is a constant or spatially varying diffusivity coefficient, and
${f(\V{x})}$
is a given right-hand side.


\subsection{Spatial discretization}
\label{sec:dg-method:discretization}

We start to describe the spatial discretization on a single-level mesh $\Omega_h$, which is obtained by decomposing the domain $\Omega$ into $\n{e}$ conforming hexahedral elements.
Further, we introduce the set of faces
${\Gamma_h}$
comprising the subsets of 
interior element interfaces $\Gamma_h^{\mathsc i}$,
Dirichlet boundary faces $\Gamma_h^{\mathsc d}$, and
Neumann boundary faces $\Gamma_h^{\mathsc n}$.
Each element $\Omega^e$ is required to possess a differentiable bijective map to the 
standard element ${\Omega^{\mathsc s} = [-1,1]^3}$ such that
${\V\xi = \V\xi^e(\V x)}$ and
${\V x  = \V  x^e(\V\xi)}$ for
${\V x  \in \Omega^e}$ and
${\V\xi \in \Omega^{\mathsc s}}$.
Given the one-dimensional polynomial basis functions 
${\{\varphi_i(\xi)\}_{i=0}^p}$,
a tensor basis of degree $p$ is constructed in $\Omega^{\mathsc s}$ 
by defining
\begin{equation*}
  \phi_{ijk}(\V\xi) = \varphi_i(\xi_1) \varphi_j(\xi_2) \varphi_k(\xi_3)
  \,.
\end{equation*}
The global basis functions in $\Omega$ are obtained by using the element mapping
and zero continuation
\begin{equation*}
  \psi_{ijk}^e(\V x) = 
    \begin{cases}
       \phi_{ijk}(\V\xi^e(\V x) & \V x \in \bar\Omega^e \\
       0                        & \V x \notin \bar\Omega^e
       \,.
    \end{cases}
\end{equation*}
This yields the trial space 
${\mathbb U_h = \operatorname{span}\bigl(\psi_{ijk}^e\bigr)}$.
The approximate solution
${u_h \in \mathbb U_h}$
is defined as
\begin{equation*}
  u_h(\V x) 
  = \sum_{e=1}^{\n{e}} \sum_{i,j,k = 0}^p u^e_{ijk} \psi_{ijk}^e(\V x)
  \qquad
  \V x \in \Omega
  \,.
\end{equation*}
This definition implies that $u_h$ is double-valued on the element interfaces,
which motivates the introduction of jump and averaging operators.
Let $\Omega^{-}$ and $\Omega^{+}$ be two adjacent elements, 
$\V n^\pm$ the outward-pointing normal vectors, and 
$u_h^\pm$ the traces of $u_h$ on ${\d\Omega^{\pm}}$.
Then the jump and the average of $u_h$ are given by
\begin{equation*}
  \jmp{u_h} = \V n^- u_h^- + \V n^+ u_h^+ 
  , \quad
  \avg{u_h} = \tfrac{1}{2}\left(u_h^- + u_h^+\right)
  \quad
  \text{on}\quad\d\Omega^{-}\!\cap\d\Omega^{+}
  \,.
\end{equation*}

While several options exist for choosing the discontinuous Galerkin method
\cite{SE_Arnold2001a}, the specific form of the method is not important 
in the following.
Here, we adopt the interior penalty method 
\cite{SE_Wheeler1978a,SE_Arnold1982a,SE_Kanschat2002a}:
Find ${u_h \in \mathbb U_h}$ such that for all ${v_h \in \mathbb U_h}$
\begin{multline}
  \label{eq:helmholtz:dgm}
  \int_{\Omega_h}
    \left( \kappa v_h u_h
         + \nu \nabla v_h \cdot \nabla u_h
    \right)
    \D\Omega 
  \\
  +
  \int_{\Gamma_h^\mathsc{i}}
    \left( \mu \max(\nu^{\pm}) \jmp{v_h}\cdot\jmp{u_h}
         - \avg{\nu\nabla v_h}\cdot\jmp{u_h}
         - \jmp{v_h}\cdot\avg{\nu\nabla u_h}
    \right) \D\Gamma
  \\
  +
  \int_{\Gamma_h^\mathsc{d}}
    \left( \mu \nu v_h 2 u_h 
         - \nu (\d_n v_h) u_h
         - \nu v_h \d_n u_h
    \right) \D\Gamma
  \\
  =
  \int_{\Omega_h}v_h f \D\Omega 
  +
  \int_{\Gamma_h^\mathsc{d}}
    \left( \mu \nu v_h 2 u_{\mathsc d} 
         - \nu (\d_n v_h) u_{\mathsc d}
    \right) \D\Gamma
  +
  \int_{\Gamma_h^\mathsc{n}}
         v_h q_{\mathsc n} \D\Gamma
  \,,
\end{multline}
where
\begin{equation*}
  \mu = \mu_{\star} \avg{\frac{p(p+1)}{2 h_{\perp}}}
\end{equation*}
is the penalty coefficient, 
${\mu_{\star} > 1}$ a non-dimensional scaling factor and
$h_{\perp}$ the element extension perpendicular to the interface
\cite{MG_Stiller2016b}.
Testing against the global basis functions yields the algebraic system
\begin{equation}
  \label{eq:helmholtz:system}
  \NM A \NM u = \NM f
\end{equation}
with
the system matrix $\NM A$,
the vector of solution coefficients ${\NM u = \bigl[u^e_{ijk}\bigr]}$ and
the right-hand side ${\NM f}$, which includes the inhomogeneous boundary conditions.

\begin{remark}[Basis functions, geometry and numerical integration]
To ease the application of overlapping Schwarz methods, a nodal basis is 
constructed using ${p\!+\!1}$ Gauss-Legendre-Lobatto points in each direction.
This basis serves to represent the approximate solution as well as the element geometry.
The corresponding quadrature is used for evaluating the element volume and surface 
integrals, which yields convergence of order ${p\!+\!1}$ with constant diffusivity 
and rectangular elements, and order $p$ in general, if the exact solution and
coefficients are smooth enough \cite{SE_Maday1990a}.
\end{remark}



\section{Multilevel formulation}
\label{sec:ml-formulation}


\subsection{Multilevel mesh}
\label{sec:ml-formulation:ml-mesh}

\begin{figure}
  \includegraphics[width=0.75\textwidth]{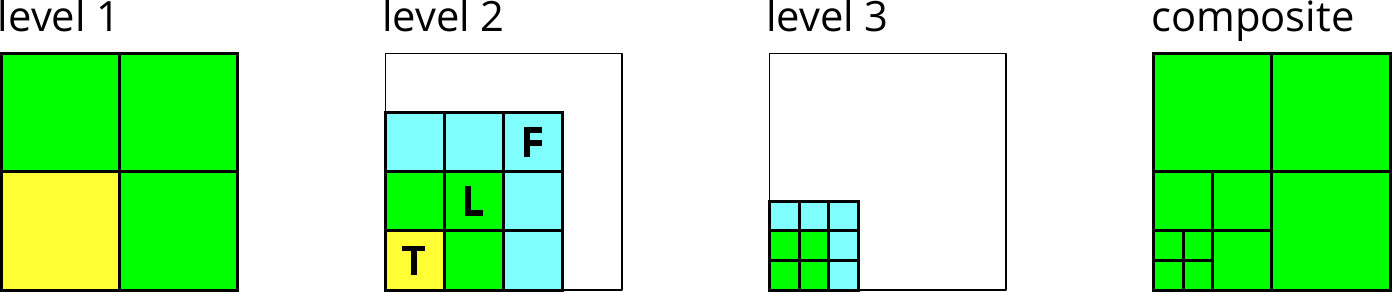}
  \caption{Locally refined multilevel mesh and corresponding composite mesh.
    Green color indicates leaf elements (L), yellow twigs (T) and cyan frozen
    elements (F). 
    \label{fig:ml-formulation:ml-mesh}}
\end{figure}

The multilevel mesh is the backbone of the multigrid method.
It is defined as a sequence of conforming hexahedral meshes 
$\{\Omega_l\}_{l=1}^{L}$
generated by successive $h$- or $p$-refinements.
With $h$-refinement, all or a subset of the elements of $\Omega_l$ are regularly subdivided to obtain the next finer mesh $\Omega_{l+1}$. 
Similarly, with $p$-refinement, the polynomial degree is elevated in the selected 
elements.
Both refinement strategies can be applied in any order, but not simultaneously.
As a result, all elements of level $l$ are assigned the same polynomial degree $p_l$.
This is primarily a technical limitation intended to preserve simple approximations and data structures.
When refining an element, two cases are distinguished:
\begin{description}
\item[regular refinement]
  yields \emph{active} child elements that participate in the solution process 
  and are coupled to the parent only via transfer operators
\smallskip
\item[irregular refinement]
  yields \emph{frozen} child elements with values  
  interpolated from parent level
\end{description}
Active elements are further distinguished into 
\emph{twigs} that are regularly refined and
\emph{leaves} that have no active children.
With global refinement, 
all elements on levels ${l < L}$ are twigs and
all elements in level ${l = L}$ are leaves.
In the case of local refinement, we enforce the following rules to maintain consistency:
\begin{enumerate}
\item
  Active elements can be marked for removal, retention or regular refinement.
\item
  Marks of elements with children selected for refinement or retention are upgraded 
  to refinement.
\item
  Neighbors of elements selected for regular refinement are marked for 
  irregular refinement, creating a layer of frozen elements that encloses 
  the active children.
\end{enumerate}
Figure~\ref{fig:ml-formulation:ml-mesh} shows a simple locally refined multilevel mesh.
Note that the root level ${l = 1}$ covers the whole domain with active elements.
Intermediate levels ${1 < l < L}$ consist of refinement zones with twigs and leaves surrounded by a layer of frozen elements.
Finally, the top level ${l = L}$ contains only leaves and frozen elements.

The function space $\mathbb U_l$, 
the trial solution $u_l$ and 
the coefficient vectors $\NM u_l$
are defined on each level as for a single mesh.
However, for the multilevel DG formulation it is essential 
to distinguish between the various classes of elements.
To this end, 
the superscripts 
$\mathsc a$,
$\mathsc f$,
$\mathsc t$ and
$\mathsc l$
are used
to refer to active, frozen, twig and leaf elements or related entities.
For example, 
${\Omega_l^\mathsc{a} = \Omega_l^\mathsc{t} \cup \Omega_l^\mathsc{l}}$
is the set of active elements on level $l$
and $\NM u_l^\mathsc{a}$ are the corresponding coefficients.
Similarly, operators are divided into suboperators, like
\begin{equation*}
\NM A_l
  = \begin{bmatrix} 
      \NM A_l^{\mathsc{a,a}} & \NM A_l^{\mathsc{a,f}} 
    \end{bmatrix}
  = \begin{bmatrix} 
      \NM A_l^{\mathsc{t,t}} & \NM A_l^{\mathsc{t,l}}  & \NM A_l^{\mathsc{t,f}} 
      \\[\smallskipamount]
      \NM A_l^{\mathsc{l,t}} & \NM A_l^{\mathsc{l,l}}  & \NM A_l^{\mathsc{l,f}}
    \end{bmatrix}
  \,,
\end{equation*}
where the first superscript indicates the result class and the second indicates 
the operand class of the respective suboperator.
Multilevel coefficients are specified by indicating the level range, 
such as
${\NM u_{l:m}}$
or simply
${\NM u_{:}}$
if all levels are included.

\subsection{Transfer operators}
\label{sec:ml-formulation:transfer-ops}

The multilevel formulation requires various operators for transferring information from coarser to finer levels and vice versa.
These operators are defined as matrices that are applied to the corresponding coefficient vectors.
Solution values or corrections given on level ${l-1}$ are prolongated to the next finer level $l$ using the embedded interpolation, which yields the coarse-to-fine interpolation operator
\begin{equation*}
  \NM u_l = \NM I_{l}\,\NM u_{l-1}
\end{equation*}
Note that the interpolation is exact due to the nested structure of the function spaces.

\begin{remark}[Notation of transfer operators]
For simplicity, the transfer operators are assigned only a single subscript that
indicates the level of result. The operand's level is omitted, as it results
from the latter and the transfer direction.
\end{remark}

The fine-to-coarse restriction operator serves for transferring residuals from the child to the current level.
As only active child elements contribute to the residual and hence only twigs receive result, the operator has the form
\begin{equation*}
  \NM r_l^\mathsc{t} = \NM R_{l}^{\mathsc{t,a}}\,\NM r_{l+1}^\mathsc{a} 
  \,,
\end{equation*}
where
$\NM R_{l}^{\mathsc{t,a}}$ 
applies to 
operands from $\Omega_{l+1}^{\mathsc a}$ and 
returns results in $\Omega_{l}^{\mathsc t}$.
Following \cite[Ch. 3.6]{MG_Hackbusch1985a}, the restriction operator is defined as the transpose of the coarse-to-fine interpolation.
Additionally, the full approximation storage multigrid scheme requires a fine-to-coarse operator for projecting solution values from child to parent levels.
Unfortunately, there is little information in the literature on how this 
operator should be defined.
In his seminal paper, \citet{MG_Brandt1977a} made no distinction between restriction and projection operators.
Later, however, he noted that these ``need not to be identical'' \cite{MG_Brandt1993a}.
\citet{MG_Trottenberg2000a} recommend the use of injection for finite difference methods and vertex-centered finite volume methods.
For high-order discontinuous Galerkin methods, embedded interpolation and $L^2$ projection proved to be equally suitable \cite{TI_Pfister2026a}.
Therefore, both options are examined here.
The resulting projection operator writes
\begin{equation*}
  \NM u_l^{\mathsc{t}} = \NM P_{l}^{\mathsc{t,a}}\,\NM u_{l+1}^\mathsc{a}
  \,.
\end{equation*}
Similar to residual restriction, only active children contribute and twigs receive a result.
In the case of $h$-refinement, eight child elements project their solution to a single parent.
To mitigate the effects of jumps, the results are averaged at shared target nodes when embedded interpolation with an even polynomial degree $p$ is used.
With odd $p$ or $L^2$ projection, no further action is required.
Alternatively, the jumps can be eliminated by subtracting a trilinear function from each child element. 
However, this approach was abandoned because it did not yield any advantage in preliminary studies.
%


\subsection{Full approximation storage formulation}
\label{sec:ml-formulation:fas}

The elliptic problem \eqref{eq:helmholtz:pde} is discretized on each level as in the single-mesh case, except for the modification of sources in twig elements and the coupling to the parent level at refinement boundaries.
In contrast to the original full approximation storage (FAS) formulation of
\citet{MG_Brandt1977a} and later work 
\cite{MG_Vanka1986a,MG_Thompson1989a,MG_Jouhaud2005a,MG_Lee2007a,MG_Feng2018a},
the coupling is not based on interpolation to the refinement boundaries.
Instead, the frozen elements are directly integrated into the DG formulation.
This yields the FAS-DG formulation 
\begin{subequations}
  \label{eq:helmholtz:fas}
  \begin{alignat}{2}
    \label{eq:helmholtz:fas:frozen}
           \NM{\tilde{u}}_l^{\mathsc{f}} 
      &  = \NM I_{l}^{\mathsc{f,l}}\,\NM{\tilde{u}}_{l-1}^{\mathsc l}
        \qquad
      && \text{in} \; \Omega^{\mathsc f}
    \\[5pt]
    \label{eq:helmholtz:fas:system}
           \NM A_l^{\mathsc{a,a}} \,\NM{\tilde{u}}_l^{\mathsc{a}}
         + \NM A_l^{\mathsc{a,f}} \,\NM{\tilde{u}}_l^{\mathsc{f}} 
      &  = \NM{\tilde{g}}_l^{\mathsc{a}}
         \qquad
      && \text{in} \; \Omega^{\mathsc a}
    \\[5pt]
    \label{eq:helmholtz:fas:rhs}
         \NM{\tilde{g}}_l^{\mathsc a} 
      & = \left\{
          \begin{array}{l}
             \smash{%
                \NM A_l^{\mathsc{t,t}} \NM P_l^{\mathsc{t,a}} 
                  \,\NM{\tilde u}_{l+1}^\mathsc{a}
              + \NM A_l^{\mathsc{t,l}} \,\NM{\tilde{u}}_l^{\mathsc{l}}
              + \NM A_l^{\mathsc{t,f}} \,\NM{\tilde{u}}_l^{\mathsc{f}}
              + \NM R_{l}^{\mathsc{t,a}}\,\NM{\tilde{r}}_{l+1}^\mathsc{a}
              }
            \\[5pt]
               \smash{\NM f_l^{\mathsc{l}}}
          \end{array}
          \right.
          \quad
      &&  \begin{array}{@{}c@{}}
             \text{in} \; \Omega^{\mathsc t} \\[5pt]
             \text{in} \; \Omega^{\mathsc l}
          \end{array}
    \\[5pt]
    \label{eq:helmholtz:fas:residual}
           \NM{\tilde{r}}_l^\mathsc{a}
      &  = \NM{\tilde{g}}_l^{\mathsc{a}}
         - \NM A_l^{\mathsc{a,af}} \,\NM{\tilde{u}}_l^{\mathsc{af}}
         \quad
      && \text{in} \; \Omega^{\mathsc a}
  \end{alignat}
\end{subequations}
where 
$\NM{\tilde{u}}_{l}$      
is the approximate numerical solution on level $l$,
$\NM{\tilde{g}}_l$
the modified RHS 
and
$\NM{\tilde{r}}_l$
the residual.
In particular, Equation \eqref{eq:helmholtz:fas:frozen} 
represents the interpolation from parent leaves to the frozen elements,
while
\eqref{eq:helmholtz:fas:system} links the active coefficients 
${\NM{\tilde{u}}_l^{\mathsc a} 
   = \transpose{[ \NM{\tilde{u}}_l^{\mathsc t} \;
                  \NM{\tilde{u}}_l^{\mathsc l} ]}
}$
with
$\NM{\tilde{u}}_l^{\mathsc{f}}$
and the right-hand side
$\NM{\tilde{g}}_l^{\mathsc{a}}$.
The composition of the right-hand side depends on the element class: 
In twigs, the operator is applied to leaf, frozen and projected child
coefficients and the result combined with the restricted child residual. 
In leaves it replicates the original RHS \eqref{eq:helmholtz:fas:rhs}.
Finally, Equation \eqref{eq:helmholtz:fas:residual} defines the residual 
as the difference between the RHS and the operator applied to the complete 
coefficient vector.

Substituting the RHS 
\eqref{eq:helmholtz:fas:rhs} in 
\eqref{eq:helmholtz:fas:system} 
yields the corresponding equations for twigs and leaves
\begin{subequations}
  \label{eq:helmholtz:fas:simplified}
  \begin{align}
    \label{eq:helmholtz:fas:simplified:twigs}
          \NM A_l^{\mathsc{t,t}} \,\NM{\tilde{u}}_l^{\mathsc{t}}
      & = \NM A_l^{\mathsc{t,t}} 
              \NM P_{l}^{\mathsc{t,a}} \,\NM{\tilde{u}}_{l+1}^\mathsc{a}
        + \NM R_{l}^{\mathsc{t,a}}\,\NM{\tilde{r}}_{l+1}^\mathsc{a}
    \\
    \label{eq:helmholtz:fas:simplified:leaves}
          \NM A_l^{\mathsc{l,t}} \,\NM{\tilde{u}}_l^{\mathsc{t}}
        + \NM A_l^{\mathsc{l,l}} \,\NM{\tilde{u}}_l^{\mathsc{l}}
        + \NM A_l^{\mathsc{l,f}} \,\NM{\tilde{u}}_l^{\mathsc{f}}
      & = \NM f_l^{\mathsc{l}}
  \end{align}
\end{subequations}
Equation~\eqref{eq:helmholtz:fas:simplified:twigs} reveals that the twig coefficients do not depend directly on the leaves and frozen elements.
This property can be exploited to reduce the cost of iterative methods by solving first \eqref{eq:helmholtz:fas:simplified:twigs} and 
then \eqref{eq:helmholtz:fas:simplified:leaves}, 
rather than tackling the system as a whole.
Upon convergence the residual vanishes and the FAS formulation simplifies to 
\begin{subequations}
  \label{eq:helmholtz:fas:converged}
  \begin{align}
	  \label{eq:helmholtz:fas:converged:frozen}
          \NM u_l^{\mathsc{f}} 
      & = \NM I_{l}^{\mathsc{f,l}} \,\NM u_{l-1}^{\mathsc l}
    \\
	  \label{eq:helmholtz:fas:converged:twigs}
          \NM u_l^{\mathsc{t}}
      & = \NM P_{l}^{\mathsc{t,a}} \,\NM u_{l+1}^\mathsc{a}
    \\
	  \label{eq:helmholtz:fas:converged:leaves}
          \NM A_l^{\mathsc{l,t}} \,\NM u_l^{\mathsc{t}}
        + \NM A_l^{\mathsc{l,l}} \,\NM u_l^{\mathsc{l}}
        + \NM A_l^{\mathsc{l,f}} \,\NM u_l^{\mathsc{f}}
      & = \NM f_l^{\mathsc{l}}
      \,.
  \end{align}
\end{subequations}
Please note that 
\eqref{eq:helmholtz:fas}
as well as 
\eqref{eq:helmholtz:fas:converged}
are fully consistent with the DG formulation 
\eqref{eq:helmholtz:system}
and therefore conservative on every level.
The coupling between levels depends on the transfer operators and is generally not conservative:
Small scales may exit a refinement zone but will not enter the enclosing region unless they can be represented on the coarser mesh.
Future work will investigate whether this is an advantage or a disadvantage.
%


1


\section{Multigrid method}
\label{sec:mg-method}


In this section, the notation is simplified for better readability:
Since the multigrid method deals only with approximations to the exact numerical 
solution, the current approximation is simply denoted by $\NM u$, while the tilde 
indicates a preliminary state, for instance, $\NM{\tilde{u}}$.


\subsection{Correction scheme}
\label{sec:mg-method:cs}

The linear multigrid method or correction scheme (CS-MG) solves the complete problem on the top level, while only corrections are computed on the coarser 
levels \cite{MG_Hackbusch1985a}.
This method has been widely used for solving linear systems arising from the discretization of elliptic equations.
In combination with global coarsening or local smoothing, CS-MG can be used also with local mesh refinement \cite{MG_Becker2000a,MG_Kanschat2004a,MG_Vu2022a,MG_Munch2023a}.
Here, the correction scheme serves as a reference for the case of global mesh refinement.
Algorithm~\ref{alg:mg-method:cs:v-cycle} outlines the V-cycle as the primary building block.
In addition to the transfer operators described in Section~\ref{sec:ml-formulation:transfer-ops}, the scheme requires methods for pre- and post-smoothing on levels ${l > 1}$
and for solving the coarse mesh problem on the root level ${l = 1}$.
The smoothers used in this study are described in Section~\ref{sec:mg-method:smoothers}.
For flexibility, the number of pre- and post-smoothing steps, $\n[l]{s1,}$ and $\n[l]{s2,}$,  can be different and vary from level to level.
The CS-MG solver consists of a sequence of V-cycles that terminates as soon as the desired accuracy (residual norm) is achieved.
Alternatively, the CS V-cycle is used as a preconditioner in an inexact preconditioned conjugate gradient method.
This approach is named CS-CG and described in \cite{MG_Stiller2016b}.

\begin{algorithm}[ht]
  \caption{Multigrid V-cycle based on the correction scheme (CS)}
  \label{alg:mg-method:cs:v-cycle}
  \begin{algorithmic}[1]
    \Procedure{CS\_V\_Cycle}{$u_L$,$f_L$}%

      \For{$l = L,2$ \,\textbf{step}\, $-1$}

        \State
          \makebox[1.8em][l]{$\NM u_l$}
            ${\; \gets \;}$ 
            ${\textsc{Smoother}(\NM{u}_l, \NM{f}_l, \n[l]{s1,})}$
        \Comment{pre-smoothing}

        \State
          \makebox[1.8em][l]{$\NM u_{l-1}$}
            ${\; \gets \;}$ 
            ${0}$
        \Comment{correction initialization}

        \State
          \makebox[1.8em][l]{$\NM f_{l-1}$}
            ${\; \gets \;}$ 
            ${\NM R_{l-1} (\NM f_l - \NM A_l \NM u_l)}$
        \Comment{residual restriction}

      \EndFor

      \State
        ${\NM u_1}$
          ${\; \gets \;}$ 
          ${\textsc{Solver}(\NM{u}_1, \NM{f}_1)}$
        \Comment{coarse mesh solution}

      \For{$l = 2,L$}

        \State
          $\NM u_l$
            ${\; \gets \;}$ 
            ${\NM u_l + \NM I_{l} \NM u_{l-1}}$
        \Comment{correction prolongation}

        \State
          $\NM u_l$
            ${\; \gets \;}$ 
            ${\textsc{Smoother}(\NM{u}_l, \NM{f}_l, \n[l]{s2,})}$
        \Comment{post-smoothing}

      \EndFor

    \EndProcedure
  \end{algorithmic}
\end{algorithm}


\subsection{Full approximation storage scheme}
\label{sec:mg-method:fas}

The full approximation scheme (FAS-MG) was developed by Achi Brandt as the backbone of his Multi-Level Adaptive Technique \cite{MG_Brandt1977a}.
As the name suggests, the method keeps an approximation to the full solution on each level.
Based on the FAS formulation (Section~\ref{sec:ml-formulation:fas}) it provides a unified approach to local mesh refinement and multigrid solution techniques.
Nevertheless, FAS-MG follows a solution strategy similar to CS-MG.
Algorithm~\ref{alg:mg-method:fas:v-cycle} presents the V-cycle adopting the FAS formulation given in Equation \eqref{eq:helmholtz:fas}.
Unlike the correction scheme, the individual steps are applied only to the involved element classes.
Note that the projection in downward leg (line 6) and the update of frozen elements in the upward leg (line 11) have no counterpart in the CS V-cycle.
Also, the coarse correction is explicitly computed before prolongating it to the present level (line 12).
While FAS-MG can be started using a full V-cycle, it is generally more efficient to employ a specific startup method.
One option is the cascadic multigrid method \cite{MG_Bornemann1996a}, which leads to Algorithm~\ref{alg:mg-method:fas:cascade}.
This method is named FAS-CMG.
It computes an initial solution on the coarsest level and interpolates it to the next finer level.
Here, the solution is improved by performing $\n[l]{sc,}$ smoothing steps and then interpolated to the next level.
This procedure is repeated until reaching the top level. 
Another option is to use the full multigrid method \cite{MG_Brandt1977a}.
The resulting method is referred to as FAS-FMG.
It starts on the root level and includes the next finer level with each iteration.
One iteration step consists of an FAS V-cycle followed by interpolation of the solution to the next level (see Algorithm~\ref{alg:mg-method:fas:fmg}).

\begin{algorithm}[t]
  \caption{Multigrid V-cycle based on the full approximation scheme (FAS)}
  \label{alg:mg-method:fas:v-cycle}
  \begin{algorithmic}[1]
    \Procedure{FAS\_V\_Cycle}{}($\NM u_{:\,}, \NM f_{:\,}$)%
      \State
        ${\NM g_{1:L} \; \gets \; \NM f_{1:L}}$

      \For{$l = L,2$ \,\textbf{step}\, $-1$}

        \State
          \makebox[1.8em][l]{$\NM u_l^{\mathsc{a}}$}
            ${\; \gets \;}$ 
            ${\textsc{Smoother}(\NM{u}_l, \NM{g}_l, \n[l]{s1,})}$
        \Comment{pre-smoothing}

        \State
          \makebox[1.8em][l]{$\NM r_l^{\mathsc{a}}$}
            ${\; \gets \;}$ 
            ${\NM g_l^{\mathsc a } 
               - \NM A_l^{\mathsc{a,a}}\,\NM u_l^{\mathsc a}
               - \NM A_l^{\mathsc{a,f}}\,\NM u_l^{\mathsc f}}$
        \Comment{residual computation}

        \State
          \makebox[1.8em][l]{$\NM u_{l-1}^\mathsc{t}$}
            ${\; \gets \;}$ 
            ${\NM P_{l-1}^{\mathsc{t,a}}\,\NM u_{l}^\mathsc{a}}$
        \Comment{solution projection}

        \State
          \makebox[1.8em][l]{$\NM g_{l-1}^\mathsc{t}$}
            ${\; \gets \;}$ 
            ${ \NM A_{l-1}^{\mathsc{t,a}}\,\NM u_{l-1}^\mathsc{a}
             + \NM A_{l-1}^{\mathsc{t,f}}\,\NM u_{l-1}^\mathsc{f}
             + \NM R_{l-1}^{\mathsc{t,a}}\,\NM r_{l}^\mathsc{a}}$
        \Comment{parent RHS update}

      \EndFor

      \State
        ${\NM u_1}$
          ${\; \gets \;}$ 
          ${\textsc{Solver}(\NM{u}_1, \NM{g}_1)}$
        \Comment{coarse mesh solution}

      \For{$l = 2,L$}

        \State
          $\NM u_l^\mathsc{f}$
            ${\; \gets \;}$ 
            ${\NM I_{l}^\mathsc{f,l}\,\NM u_{l-1}^\mathsc{l}}$
        \Comment{frozen element update}

        \State
          $\NM u_l^\mathsc{a}$
            ${\; \gets \;}$ 
            ${ \NM u_l^\mathsc{a} 
             + \NM I_{l}^\mathsc{a,t}
                 ( \NM u_{l-1}^\mathsc{t} 
                 - \NM P_{l-1}^\mathsc{t,a}\,\NM u_l^\mathsc{a} )
            }$
        \Comment{correction prolongation}

        \State
          $\NM u_l^\mathsc{a}$
            ${\; \gets \;}$ 
            ${\textsc{Smoother}(\NM{u}_l, \NM{g}_l, \n[l]{s2,})}$
        \Comment{post-smoothing}

      \EndFor

    \EndProcedure
  \end{algorithmic}
\end{algorithm}

\begin{algorithm}[ht]
  \caption{Cascadic multigrid method (CMG)}
  \label{alg:mg-method:fas:cascade}
  \begin{algorithmic}[1]
    \Procedure{FAS\_Cascade}{}($\NM u_{:\,}, \NM f_{:\,}$)%

      \State
        ${\NM u_1}$
          ${\; \gets \;}$ 
          ${\textsc{Solver}(\NM{u}_1, \NM{f}_1)}$
        \Comment{coarse mesh solution}

      \For{$l = 2,L-1$}

        \State
          $\NM u_l$ ${\; \gets \;}$ ${\NM I_{l}\,\NM u_{l-1}}$
        \Comment{interpolation}

        \State
          $\NM u_l^\mathsc{a}$
            ${\; \gets \;}$ 
            ${\textsc{Smoother}(\NM{u}_l, \NM{f}_l, \n[l]{sc,})}$
        \Comment{approximate solution}

      \EndFor

    \EndProcedure
  \end{algorithmic}
\end{algorithm}

\begin{algorithm}[ht]
  \caption{Full multigrid method (FMG)}
  \label{alg:mg-method:fas:fmg}
  \begin{algorithmic}[1]
    \Procedure{FAS\_FMG}{}($\NM u_{:\,}, \NM f_{:\,}$)%

      \For{$l = 1,L-1$}

        \State
          \makebox[1.8em][l]{$\NM u_{1:l}$}
            ${\; \gets \;}$ \textsc{FAS\_V\_Cycle}($\NM u_{1:l}, \NM f_{1:l}$)
     
        \State
          \makebox[1.8em][l]{$\NM u_{l+1}$} 
            ${\; \gets \;}$ ${\NM I_{l+1}\,\NM u_{l}}$

      \EndFor

    \EndProcedure
  \end{algorithmic}
\end{algorithm}


\subsection{Smoothers and solvers}
\label{sec:mg-method:smoothers}

During a multigrid cycle, linear systems of the form
\begin{equation}
  \label{eq:mg-method:smoothers:system:solution}
  \NM A \NM u = \NM g
\end{equation}
must be smoothed or solved repeatedly on every level.
Note that level subscripts and element classifiers have been omitted for readability.
Starting from the initial approximation $\NM{\tilde u}$, 
Equation \eqref{eq:mg-method:smoothers:system:solution} 
can be transformed into a linear system for the correction 
${\Delta \NM u = \NM u - \NM{\tilde u}}$
with the residual
${\NM r = \NM g - \NM A \NM{\tilde u}}$
serving as the RHS:
\begin{equation}
  \label{eq:mg-method:smoothers:system:correction}
  \NM A \Delta \NM u = \NM r
  \,.
\end{equation}
On the finer levels, ${l > 1}$, a smoothing method is applied to
Equation~\ref{eq:mg-method:smoothers:system:correction}.
Its main purpose is to eliminate the short-wave error components from 
$\NM{\tilde u}$, while the remaining long-wave components are removed
recursively on the coarser meshes \cite{MG_Hackbusch1985a,MG_Bramble1993a}.
The following subsections describe the Schwarz smoother and the IPCG method, 
which is used both as a smoother and as a solver on the root level.
%


\subsubsection{Weighted additive Schwarz method}
\label{sec:mg-method:smoothers:schwarz}

Schwarz methods are iterative substructuring methods based on domain decomposition \cite{MG_Dryja1994a}.
In the spectral element community, these methods are often combined with fast diagonalization (FD) and multigrid methods.
The FD method goes back to \citet{MG_Lynch1964a} and was adapted to SEM by \citet{MG_Couzy1994a}, while
\citet{MG_Pahl1993a} and
\citet{MG_Casarin1997a}
pioneered overlapping Schwarz methods.
Combining these ideas, \citet{MG_Fischer2000a} presented
a two-level Schwarz preconditioner using 
the FD method for solving the subproblems and 
a direct coarse grid solver for eliminating long-wave residual components.
Using the Schwarz method as a smoother in multigrid methods, 
\citet{MG_Lottes2005a} developed iterative solvers that achieve
superior efficiency and robustness.
\citet{MG_Stiller2016a} achieved a further improvement by introducing
a non-uniform weighting of the subdomain solutions, reducing the number 
of multigrid cycles by a factor of up to 3.
Nearly identical results were obtained with discontinuous Galerkin methods
using rectangular elements \cite{MG_Stiller2016b,MG_Stiller2017a}.

\begin{figure}
  \subcaptionbox{subdomain $\Omega^e_{\mathsc s}$ on the unstructured curved mesh
    \label{fig:mg-method:smoothers:schwarz:orig}}
    {\includegraphics[scale=0.55]{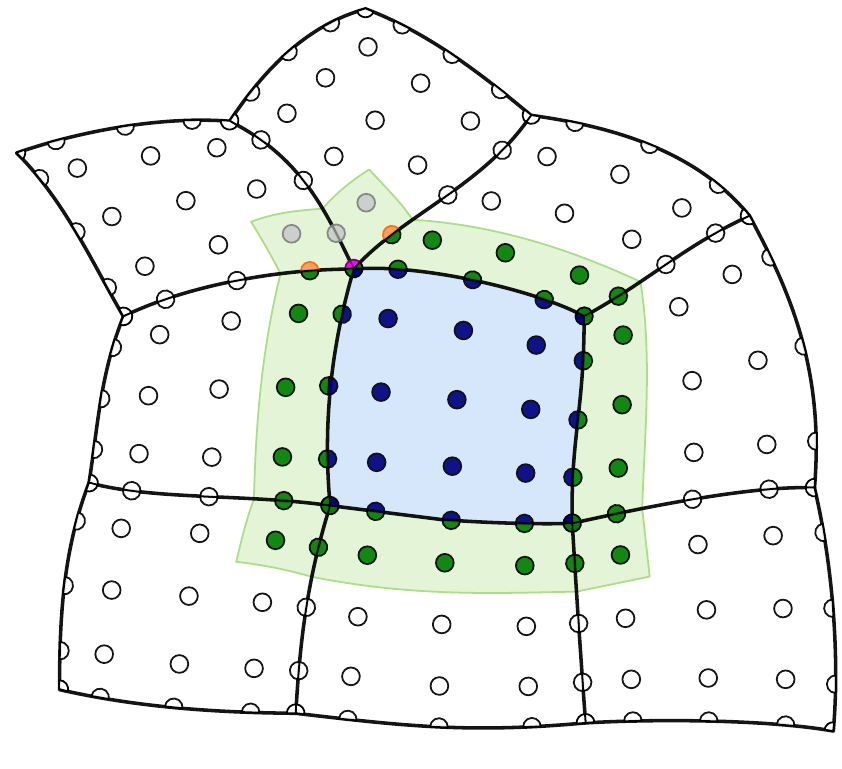}}
  \hfill
  \subcaptionbox{Cartesian surrogate domain $\tilde\Omega^e_{\mathsc s}$
    \label{fig:mg-method:smoothers:schwarz:approx}}
    {\includegraphics[scale=0.55]{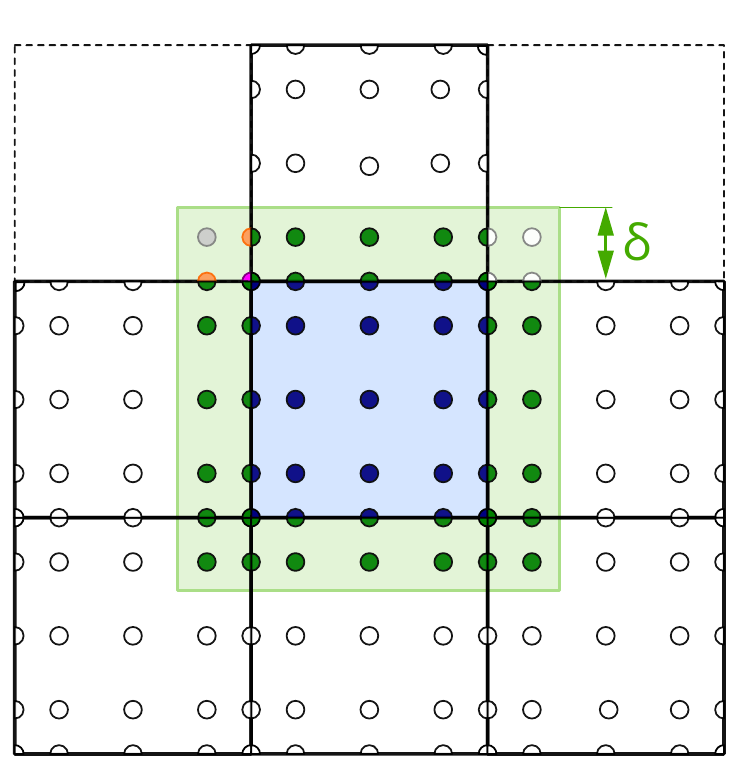}}
  \caption{Two-dimensional illustration of an element centered overlapping subdomain. 
    \label{fig:mg-method:smoothers:schwarz:subdomain}}
\end{figure}

Only a few studies addressed the generalization to unstructured curvilinear meshes.
For the case of deformed geometries \citet{MG_Couzy1994a} proposed to approximate 
element subproblems on the corresponding rectangular elements so that the FD method 
remains applicable.
The resulting method is no longer an exact solver to the subproblems, but still a good preconditioner for the overall problem.
\citet{MG_Fischer2005a} extended this approach to overlapping Schwarz methods for continuous spectral elements on unstructured meshes.
To cope with irregular mesh topologies, they constructed the subdomain problems using only the neighbors sharing a common face.
Entries corresponding to overlapped edge and vertex neighbors were set to zero.
It should be noted, however, that these neighbors remain connected via shared nodes.
\citet{MG_Vu2022a} 
developed an overlapping Schwarz preconditioner for DG methods that extracts 
the subdomains directly from the curvilinear mesh.
While this approach promises higher robustness against deformations, it prevents
the use of the FD method.
Instead, the authors apply an ILUT method, which is still efficient but more costly 
than FD and requires a specific preprocessing for each subdomain.
To reduce complexity, the method includes only face neighbors in the subdomains.
However, this simplification can compromise robustness, as the DG method lacks a direct coupling to edge and vertex neighbors.
We adopt a different approach that generalizes the concepts put forward in
\cite{ %
  MG_Couzy1994a,MG_Fischer2000a,MG_Lottes2005a, %
  MG_Stiller2016a,MG_Stiller2016b,MG_Stiller2017a}.
The design principles of the proposed Schwarz method are as follows:
\begin{enumerate}
\item
  The computational domain is partitioned into overlapping, element-centered 
  subdomains $\Omega^e_{\mathsc s}$, which have a one-to-one correspondence to
  the spectral elements $\Omega^e$.
\item
  Each subdomain consists of the central element and an overlap region including 
  adjoining parts of the surrounding elements (blue or green shaded areas in
  Figure~\ref{fig:mg-method:smoothers:schwarz:orig}).
\item
  The overlap is determined by the parameter ${0 \le \delta < 1}$,
  which specifies the normalized width in the standard element and thus defines
  the number of node layers to be included
  (Figure~\ref{fig:mg-method:smoothers:schwarz:approx}).
\item
  The corresponding layers are adopted from all neighbors for which there exists
  a one-to-one mapping between the original and the approximate subdomain
  (green nodes in Figure~\ref{fig:mg-method:smoothers:schwarz:subdomain}).
\item
  Edge or vertex neighbors that have no counterpart in the original subdomain 
  are filled with zero entries (symbolized by empty nodes in the upper right 
  corner of the overlap region in Figure~\ref{fig:mg-method:smoothers:schwarz:approx}).
\item
  If there exist more neighbors than in the regular case, the innermost edge or vertex 
  nodes are arithmetically weighted and merged into a single node in the approximate
  subdomain (magenta nodes at the upper left corner of the central element in 
  Figure~\ref{fig:mg-method:smoothers:schwarz:subdomain}).
\item
  Optionally, irregular edge nodes adjoining face neighbors can be included
  (colored orange in Figure~\ref{fig:mg-method:smoothers:schwarz:subdomain}).
\item
  The original subproblem is approximated using Cartesian tensor-product operators
  with averaged coefficients.
\item
  The resulting regular problems are solved using the FD method.
\item
  The subdomain solutions are weighted and combined to global corrections
  as described in \cite{MG_Stiller2016a,MG_Stiller2016b}.
\end{enumerate}
Please note that, unlike the approach taken by \citet{MG_Fischer2000a}, the regular edge and vertex neighbors are retained and treated in the same way as in the structured case.

Following the design principles, the subproblem in $\Omega^e_{\mathsc s}$
is defined as
\begin{equation}
  \label{eq:mg-method:smoothers:schwarz:exact}
  \NM A^e_{\mathsc s} \Delta \NM u^e_{\mathsc s}
    = \NM r^e_{\mathsc s}
    \,
\end{equation}
where
$\NM R^e_{\mathsc s}$
is the operator for node restriction and merging 
(not to be confused with the fine-to-coarse restriction operator),
${
  \NM A^e_{\mathsc s} 
    = \NM R^e_{\mathsc s}
      \NM A
      \transpose{(\NM R^e_{\mathsc s})}
}$
the subdomain operator,
$\NM r^e_{\mathsc s}$ 
the subdomain residual and
$\Delta \NM u^e_{\mathsc s}$
the subdomain correction.
To apply the FD method, the subdomain $\Omega^e_{\mathsc s}$ on the unstructured curved mesh 
is approximated by a rectangular surrogate domain $\tilde\Omega^e_{\mathsc s}$, which is obtained in four steps:
\begin{enumerate}
\item
  transform the isoparametric map of the central element to Legendre polynomials, 
\item
  drop all nonlinear terms to obtain the approximate parallelepiped,
\item
  define the corresponding cuboid by neglecting non-orthogonality and 
  rescaling the edge lengths to ${\Delta x^e_i}$ so that the element volume 
  is preserved,
\item
  construct a ${3 \times 3 \times 3}$ Cartesian mesh by attaching further cuboids
  as depicted in Figure~\ref{fig:mg-method:smoothers:schwarz:approx}.
\end{enumerate}
In addition, variable diffusion coefficients are replaced by their element average 
$\tilde{\nu}^e$.
Based on these simplifications, the approximate subdomain operator is defined as
\begin{equation}
  \label{eq:mg-method:smoothers:schwarz:tpo}
  \begin{alignedat}{2}
    \NM{\tilde{A}}^e_{\mathsc s}
      &  = g_0\,
           \NM M_{3} \otimes
           \NM M_{2} \otimes
           \NM M_{1}
      && + g_1\,
           \NM M_{3} \otimes
           \NM M_{2} \otimes
           \NM L_{1}
      \\
      &  + g_2\,
           \NM M_{3} \otimes
           \NM L_{2} \otimes
           \NM M_{1}
      && + g_3\,
           \NM L_{3} \otimes
           \NM M_{2} \otimes
           \NM M_{1}
           \,,
  \end{alignedat}
\end{equation}
where
\begin{equation*}
  g_0 = \kappa \frac{\Delta x^e_{1}\Delta x^e_{2}\Delta x^e_{3}}{8},
  \quad
  g_i = \tilde{\nu}^e \frac{\Delta x^e_{j}\Delta x^e_{k}}{2\Delta x^e_{i}}
  \quad
  (i,j,k \ \text{cyclic}).
\end{equation*}
The mass matrices 
$\NM M_i$
are symmetric positive-definite
and stiffness matrices
$\NM L_i$
symmetric positive-semidefinite.
While the approximate operator 
$\NM{\tilde{A}}^e_{\mathsc s}$ 
as a whole is specific to each element, 
the mass and stiffness matrices
depend only on the conditions at the corresponding faces.
On each side of the subdomain either
interior ($\mathsc i$),
Dirichlet ($\mathsc d$) or
Neumann ($\mathsc n$)
conditions can be imposed.
Thus, for any system matrix ${\NM X \in \{\NM M_i, \NM L_i \}}$
there exists a set of nine variants such that
\begin{equation*}
  \NM X
  \in
  \{ \NM X^\mathsc{ii}
   , \NM X^\mathsc{id}
   , \NM X^\mathsc{in}
   , \NM X^\mathsc{di}
   , \NM X^\mathsc{dd}
   , \NM X^\mathsc{dn}
   , \NM X^\mathsc{ni}
   , \NM X^\mathsc{nd}
   , \NM X^\mathsc{nn}
  \}
  \,,
\end{equation*}
where the superscripts indicate the conditions on the left and the right side, respectively.
Explicit representations of the corresponding mass and stiffness matrices can be found in \cite{MG_Stiller2016b}.

Each pair of mass and stiffness matrices is associated with a 
generalized eigenvalue problem of the form
${
  \NM L\, \NM s  = \lambda \NM M \NM s 
}$,
which leads to the similarity transform
\begin{equation*}
  \transpose{\NM S} \NM L\,\NM S = \NM \Lambda, \quad
  \transpose{\NM S} \NM M  \NM S = \NMC I
\end{equation*}
with
the eigenvectors
${\NM S = [ \NM s_i]}$
and the eigenvalues
${\NM \Lambda = \diag \lambda_i}$,
see \cite{MG_Lynch1964a}.
Using these properties and the tensor-product structure of
\eqref{eq:mg-method:smoothers:schwarz:tpo},
the FD method yields the inverse operator
\begin{equation*}
     ( \NM{\tilde{A}}^e_{\mathsc s} )^{-1}
     = ( \NM S_3 \otimes \NM S_2 \otimes \NM S_1 )
       \NM D^{-1}
       ( \transpose{\NM S_3} \otimes 
         \transpose{\NM S_2} \otimes 
         \transpose{\NM S_1} )
\end{equation*}
with the diagonal matrix
\begin{equation*}
  \NM D = g_0 \, \NMC I_3 \otimes \NMC I_2 \otimes \NMC I_1
        + g_1 \, \NMC I_3 \otimes \NMC I_2 \otimes \NM \Lambda_1
        + g_2 \, \NMC I_3 \otimes \NM \Lambda_2 \otimes \NMC I_1
        + g_3 \, \NM \Lambda_3 \otimes \NMC I_2 \otimes \NMC I_1
        \,.
\end{equation*}
Note that the eigenvectors $\NM S$ and eigenvalues $\NM \Lambda$
depend only on the face conditions and, hence, can be extracted
from a small set of precomputed configurations.
Applying the inverse operator to Equation~\eqref{eq:mg-method:smoothers:schwarz:exact}
yields the approximate subdomain correction
\begin{equation*}
  \Delta \NM u^e_{\mathsc s}
    \simeq ( \NM{\tilde{A}}^e_{\mathsc s} )^{-1} \NM r^e_{\mathsc s}
    \,.
\end{equation*}
Due to the tensor-product structure, the evaluation of the correction is 
very efficient and requires only $O(p)$ operations per subdomain node.
Finally, the overall correction results as the weighted combination of all subdomain 
corrections, i.e.
\begin{equation}
  \label{eq:mg-method:smoothers:schwarz:correction}
  \Delta u = \sum_{e}^{\n{e}} 
                \transpose{(\NM R^e_{\mathsc s})} \,
                \NM W^e \Delta \NM u^e_{\mathsc s}
           = \sum_{e}^{\n{e}} 
                \transpose{(\NM R^e_{\mathsc s})} \,
                \NM W^e 
                ( \NM{\tilde{A}}^e_{\mathsc s} )^{-1} \,
                \NM R^e_{\mathsc s}\;
                \NM r
                \,,
\end{equation}
where
\begin{equation*}
  \NM W^e = \NM W_3 \otimes \NM W_2 \otimes \NM W_1
\end{equation*}
is a diagonal tensor-product weighting matrix that depends only on the conditions 
imposed at the sides of the subdomain.
To compute the weights, we use the quintic weighting function introduced in 
\cite{MG_Stiller2016a}. 
The iterative application of 
Equation~\eqref{eq:mg-method:smoothers:schwarz:correction}
leads to the weighted Schwarz method described in 
Algorithm~\ref{alg:mg-method:smoothers:schwarz}.
This method is referred to as WAS.

\begin{algorithm}[ht]
  \caption{Weighted additive Schwarz method (WAS)}
  \label{alg:mg-method:smoothers:schwarz}
  \begin{algorithmic}[1]
    \Procedure{WAS\_Method}{}($\NM u, \NM f, \n{i}$)%

      \For{$i = 1,\n{i}$}

        \State
          \makebox[0.5em][l]{$\NM r$}
            ${\; \gets \;}$ 
            ${\NM f - \NM A \NM u}$
          \Comment{residual computation: $\NM r = \transpose{[\NM r^{\mathsc{a}} \; \NMC 0]}$}

        \State
          \makebox[0.5em][l]{$\NM u$}
            ${\; \gets \;}$ 
            ${ \NM u 
             + \sum_{e=1}^{\n{e}} 
                  \transpose{(\NM R^e_{\mathsc s})} \,\NM W^e 
                  ( \NM{\tilde{A}}^e_{\mathsc s} )^{-1} \,\NM R^e_{\mathsc s}\; \NM r
            }$
        \Comment{solution update, $\NM u^{\mathsc f}$ unchanged}

      \EndFor

    \EndProcedure
  \end{algorithmic}
\end{algorithm}


\subsubsection{Inexact preconditioned conjugate gradient method}
\label{sec:mg-method:smoothers:ipcg}

The Schwarz method has proved to be an excellent smoother for multigrid methods 
used with high-order spectral element methods
\cite{MG_Lottes2005a,MG_Stiller2016b,MG_Vincent2019a}.
Unfortunately, the method is vulnerable to mesh anisotropy and irregularity,
which are common in complex applications.
With global refinement, this problem can be mitigated by Krylov 
acceleration, i.e. by using the Schwarz multigrid method as a preconditioner
for a conjugate gradient (CG) method or a GMRES method
\cite{MG_Lottes2005a,MG_Stiller2016b}.
A further option is to accelerate the Schwarz method using a Chebyshev method
\cite{MG_Phillips2025a}.
To be efficient, the Chebyshev method requires the knowledge of the largest
eigenvalue of the preconditioned system matrix.
However, as a result of mesh adaptation, the eigenvalue changes and
must be computed repeatedly.
Therefore, this approach appears to be unsuitable for adaptive multigrid methods.
Instead, we adopt the inexact preconditioned conjugate gradient method of \citet{KR_Golub1999a}.
Using the method with one Schwarz sweep for preconditioning yields the
Schwarz-preconditioned conjugate gradient (SCG) method described in
Algorithm~\ref{alg:mg-method:smoothers:scg}.
Note that the Schwarz method is inlined instead of using the procedure
defined in Algorithm~\ref{alg:mg-method:smoothers:schwarz}.
The SCG method is utilized both as a smoother and as a solver for the 
coarse problem on the root level.

\begin{remark}[IPCG for semidefinite systems]
Like the standard CG method, IPCG requires the system matrix to be positive-definite.
However, in the present case $\NM A$ is only semidefinite
if ${\kappa}$ vanishes and no Dirichlet conditions are prescribed.
Nevertheless, convergence can be ensured by projecting RHS to the null space of $\NM A$
(see \citet{KR_Kaasschieter1988a}).
\end{remark}

\begin{algorithm}[ht]
  \caption{Inexact conjugate gradient method with Schwarz preconditioner (SCG)}
  \label{alg:mg-method:smoothers:scg}
  \begin{algorithmic}[1]
    \Procedure{SCG\_Method}{}($\NM u, \NM f, \n{i}, r_{\mathsc{max}}$)%

      \State
        \makebox[0.5em][l]{$\NM r$} ${\; \gets \; \NM f - \NM A \NM u}$
        \Comment{initial residual: $\NM r^{\mathsc f} = \NMC 0$}
      \State
        \makebox[0.5em][l]{$\NM p$} ${\; \gets \; \NM r}$

      \For{$i = 1,\n{i}$}

        \State
          \makebox[0.5em][l]{$\NM z$}
            ${\; \gets \;}$ 
            ${\sum_{e=1}^{\n{e}} 
                 \transpose{(\NM R^e_{\mathsc s})} \,\NM W^e 
                 ( \NM{\tilde{A}}^e_{\mathsc s} )^{-1} \,\NM R^e_{\mathsc s}\; \NM r
            }$
          \Comment{Schwarz sweep}

        \If{$i > 1$}
          \State
            \makebox[0.5em][l]{$\beta$} 
              ${\; \gets \; \transpose{\NM z}(\NM r - \NM r_0) \, / \, \gamma}$
          \State
            \makebox[0.5em][l]{$\NM p$} ${\; \gets \; \beta\NM p + \NM z}$                
        \EndIf

        \State
          \makebox[0.8em][l]{$\NM q$} ${\; \gets \; \NM A \NM p}$
        \State
          \makebox[0.8em][l]{$\gamma$} 
            ${\; \gets \; \transpose{\NM z} \NM r}$
        \State
          \makebox[0.8em][l]{$\alpha$} 
            ${\; \gets \; \gamma \, / \, \transpose{\NM p} \NM q}$
        \State
          \makebox[0.8em][l]{$\NM r_0$} ${\; \gets \; \NM r}$
        \State
          \makebox[0.8em][l]{$\NM r$} ${\; \gets \; \NM r - \alpha \NM q}$
          \Comment{residual update, $\NM r^{\mathsc f}$ unchanged}
        \State
          \makebox[0.8em][l]{$\NM u$} ${\; \gets \; \NM u + \alpha \NM p}$
          \Comment{solution update, $\NM u^{\mathsc f}$ unchanged}
        \State
          \textbf{if} \ $\norm{r} \le r_{\mathsc{max}}$ \ \textbf{exit}

      \EndFor

    \EndProcedure
  \end{algorithmic}
\end{algorithm}


\subsection{Cost analysis}
\label{sec:mg-method:cost}

The cost of multigrid methods depends on various factors:
\begin{itemize}
\item
  the number of V-cycles $\n{v}$,
\item
  the smoother,
\item
  the number of pre-smoothing steps $\n{s1}$ and post-smoothing steps $\n{s2}$,
\item
  the number of top-level smoothing steps $\n{sc}$ between consecutive V-cycles,
\item
  the transfer operators,
\item
  the startup method (for FAS-MG), and
\item
  the cost of Krylov acceleration on the top level (for CS-CG).
\end{itemize}

The cost of a single operator application on level $l$ is estimated as
\begin{equation*}
  w_{\mathsc a,l} = c_{\mathsc a} n_l p_l
  \,,
\end{equation*}
where 
$c_{\mathsc a}$ is a constant and
${n_l = \n[l]{e,} (p_l+1)^3}$ is the number of unknowns on level $l$.
Ignoring the overlap between the subdomains, the cost of the smoother is
\begin{equation*}
  w_{\mathsc s,l} = c_{\mathsc s} n_l p_l
  \,.
\end{equation*}
For simplicity we assume 
${c_{\mathsc s} = c_{\mathsc a}}$ for WAS and
${c_{\mathsc s} = 2 c_{\mathsc a}}$ for SCG.
On the root level the cost of the solver is estimated by
${\n{s1} + \n{s2}}$ applications of the smoother.
The computation of the parent RHS entails evaluating the residual, 
which in turn requires one operator application.
Note that FAS involves an additional operator application on the parent level.
However, due to the coarser mesh, this operation is significantly less expensive and consequently neglected.
Similarly, the cost of interpolation, restriction and projection operators is ignored.
Under these assumptions, the cost of $\n{i}$ iterations can be estimated as
\begin{equation}
  \label{eq:mg-method:cost:total}
  W(c_{\mathsc s},c_{\mathsc k},\n{s},\n{sc},\n{i})
    = n_L p_L
      \Bigl[ \n{i} \bigl(\n{s} c_{\mathsc s} + c_{\mathsc a})
             \frac{1 - r_{\mathsc a}^L}{1 - r_{\mathsc a}}
           + (\n{i} - 1) (\n{sc} - \n{s})c_{\mathsc s}
           + \n{i} c_{\mathsc k}
      \Bigr],
\end{equation}
where
${\n{s} = \n[1]{s} + \n[2]{s}}$
is the number of pre- and post-smoothing steps per level,
${r_{\mathsc a} \simeq w_{\mathsc a,l} / w_{\mathsc a,l+1}}$
is the operator cost ratio, 
and 
$c_{\mathsc k}$ the cost factor for top-level Krylov acceleration, 
which is $c_{\mathsc a}$ for CS-CG and $0$ else.
Note that, in general, the operator cost ratio varies from level to level.
Here, we assume a constant coarsening ratio ${q = n_l/n_{l+1}}$, which results in
${r_{\mathsc a} = q}$ for $h$-MG and
${r_{\mathsc a} = q/2}$ for $p$-MG.
Typical values of $q$ are
$1/8$ for global refinement and
$1/4$ for local refinement near surfaces.

The cost of the startup methods can be estimated as
\begin{equation*}
  W_{\mathsc{CMG}}
    < n_L p_L \n[0]{s} c_{\mathsc s} \frac{a}{1 - a}
\end{equation*}
for cascadic multigrid and, respectively,
\begin{equation*}
  W_{\mathsc{FMG}}
    < n_L p_L (\n{s} c_{\mathsc s} + c_{\mathsc a})
      \frac{a}{(1 - a)^2}
\end{equation*}
for full multigrid.
Consequently, the startup method requires only a small fraction of the cost of a single V-cycle and is therefore neglected.

For a fair comparison between different methods, the total cost 
\eqref{eq:mg-method:cost:total} 
is normalized by the cost of a single V-cycle using the
WAS smoother with ${\n[1]{s} = \n[2]{s} = 1}$.
This yields the equivalent number of \emph{standard} V-cycles
\begin{equation}
  \label{{eq:mg-method:cost:equivalent}}
  v
    = \frac{W(c_{\mathsc s},c_{\mathsc k},\n{s},\n{sc},\n{i})}
           {W(c_{\mathsc a},0,2,2,1)}
    = \n{i} 
      \frac{\n{s} c_{\mathsc s} + c_{\mathsc a}}
           {3 c_{\mathsc a}}
    + \frac{(\n{i} - 1)(\n{sc} - \n{s})c_{\mathsc s} + \n{i}c_{\mathsc k}}
           {3 c_{\mathsc a}}
      \frac{1 - a}{1 - a^L}
      \,.
\end{equation}
%





\section{Data structures, adaptation and parallelization}
\label{sec:implementation}


\subsection{Data structures}
\label{sec:implementation:structures}

When designing data structures, it is advisable to separate 
the basic spectral element operators,
the multilevel mesh and
the problem-specific data.
The multilevel mesh comprises a hierarchy of partitioned meshes, 
which in turn consist of elements.
While the mesh forms the backbone of the method, 
the elements are the fundamental building blocks.
Since there are many ways to structure the element data type, we list here only the components that are essential for implementing the method:
\begin{itemize}
\item
   state: active or frozen
\item
   for all faces, edges and vertices: 
   list of neighbors, including
   \begin{itemize}
   \item
     element ID
   \item
     partition ID
   \item
     related face, edge or vertex
   \item
     relative orientation
   \end{itemize}
\item
   parent: element ID and partition ID
\item
   children:
   \begin{itemize}
   \item
     refinement type: regular, irregular, or none
   \item
     first child element ID and partition ID
   \end{itemize}
\end{itemize}
The neighbor list includes ghosts, which refer to elements located in other partitions.
Nevertheless the ghosts are assigned an ID beyond the range of local elements, which simplifies the handling of related mesh or solution data.
The neighbor orientation is necessary for composing interface fluxes and restricting data to overlapping subdomains.
It is always required for faces and edges, but may be omitted for irregular vertices.
The parent and child data serve for constructing transfer operators but also for 
controlling refinement and remapping during mesh adaptation.
Depending on the implementation, more information may be included such as geometry data, 
face boundary conditions or the position on a space filling curve (SFC).
Alternatively, selected element data can be stored in separate structures like hash tables.

We would also like to point out that some of the information and data management tasks 
can be outsourced to external libraries that are optimized for this purpose, such as 
p4est \cite{MG_Burstedde2011a}.


\subsection{Adaptation and partitioning}
\label{sec:implementation:partitioning}

Adaptation involves refining or removing elements to meet changing resolution 
requirements.
Typically, this results in an imbalance, which must be resolved by redistributing 
the multilevel mesh. 
Both steps present a challenge that can be met in different ways.
In this study, the mesh is adapted, rebuilt, and redistributed level by level.
Partitioning can be based on abstract weighted graphs or on space-filling curves 
\cite{MG_Karypis1998a,MG_Boman2012a,MG_Griebel1999a,MG_Burstedde2011a,
      MG_Weinzierl2019a,MM_Cerveny2019a,MG_Clevenger2020a}.
Here we use the latter approach.
For this purpose, all root level elements are ordered along a Hilbert polygon
defined in the bounding box of the computational domain.
The required mapping between physical and SFC coordinates is provided by the
FD4 library \cite{MG_Lieber2010a,MG_Lieber2018a}.
The SFC ordering is continued on subsequent levels by subdividing the polygon
within parent elements, eliminating the need for recursion or communication.
Once a consistent adaptation pattern has been established, the root level
is partitioned and redistributed first. 
Then the higher levels are rebuilt recursively step by step.
First, the elements on the present level create prototypes of their children.
Simultaneously the SFC is subdivided and used to assign the children to their
target partitions.
To simplify transfer operations, the children of each element are grouped into 
a cluster and thus directed to the same destination.
In the next step, the prototypes are sent to their destinations, where they are 
completed and combined with elements from other parents to form a new partition.
Finally, metadata -- such as adaptation marks -- and solution data are transferred 
from the retained elements to the corresponding new instances.
Regardless of its complexity, grid adaptation is less expensive than a V-cycle.


\subsection{Implementation}
\label{sec:implementation:hispeet}

The methods presented in this study have been implemented 
in the high-order spectral-element library HiSPEET.
The sources including examples is available at
\url{https://github.com/HiSPEET/hispeet.git}.
HiSPEET builds on the FD4 library
\cite{MG_Lieber2018a} 
for Hilbert curves and on LIBXSMM
\cite{HPC_Georganas2021_libxsmm}
for the fast evaluation of spectral element operators.




\section{Numerical experiments}
\label{sec:num-exp}


\subsection{Overview}
\label{sec:num-exp:notation}

This section presents numerical experiments illustrating the performance of
the FAS-MG method with high-order DG methods.
The first part focuses on optimizing smoothing methods and cycle parameters 
to ensure robust execution on anisotropic, curved, and irregular meshes.
Following this, the performance is evaluated on locally refined meshes.
The third part examines the parallel adaptive solution of a test case 
involving a thin spherical front.
Finally, we consider the application of locally refined meshes in a flow solver
as an example extending the approach to more complex problems.

The methods introduced in Section~\ref{sec:mg-method} are specified in more detail 
by appending the name of the smoother and the number of smoothing steps used.
Further, the latter are distinguished into 
$\n{s1}$ pre-smoothing steps,
$\n{s2}$ post-smoothing steps, and
$\n{sc}$ crossover-smoothing steps, applied on the top level between two consecutive cycles.
For example, FAS-FMG-WAS(3,1,2) denotes the FAS-MG method starting with full multigrid
and using the WAS smoother with
${\n{s1} = 3}$,
${\n{s2} = 1}$ and
${\n{sc} = 2}$
steps for pre-, post- and crossover smoothing.

The error ${\varepsilon_h = u_h - u}$ is measured using the $L^2$ norm 
on the composite mesh
\begin{equation*}
  \varepsilon_0 
  = \biggl(
    \,
      \sum_{l=1}^{L} \int_{\Omega_{l}^{\mathsc l}} \varepsilon_{l}^2 \D\Omega
    \,
    \biggr)^{1/2}
\end{equation*}
or the energy norm
\begin{equation*}
  \varepsilon_{\mathsc e} 
  = \biggl(
    \,
      \sum_{l=1}^{L} \,
        \transpose{\bigl(\,\NM{\varepsilon}_{l}^{\mathsc l}\,\bigr)} 
        \, \NM A_{l}^{\mathsc{l,af}} \NM{\varepsilon}_{l}
    \,
    \biggr)^{1/2}
    \,.
\end{equation*}
To residuals we apply the Euclidian norm 
\begin{equation*}
  r_2 
  = \biggl(
    \,
      \sum_{l=1}^{L} \,
        \transpose{\bigl(\,\NM{r}_{l}^{\mathsc l}\,\bigr)} 
        \, \NM{r}_{l}^{\mathsc l}
    \,
    \biggr)^{1/2}  
\end{equation*}

All numerical experiments were run on TU Dresden's HPC system Barnard. 
Each node of Barnard consists of two Intel Xeon Platinum 8470 CPUs with 52 cores
clocked 2.00 GHz. 
The code is compiled with the GNU compiler collection 12.2.0
and parallelized with Open MPI 4.1.4.
Start values are generated using the inbuilt random number generator and scaled to
the interval ${[-1,1]}$.
The performance of the multigrid methods is evaluated using
\begin{itemize}
\item
  the number $n_{10}$ of cycles or iterations to reduce the residual 
  by a factor of $10^{10}$,
\item
  the equivalent number of standard V-cycles 
  ${v_{10} = v(n_{10})}$,
\item
  the corresponding wall time $t_{10}(\n{p})$ using $\n{p}$ MPI ranks
  in seconds, and
\item
  the wall time $\tau_{10}(\n{p})$ per core and unknown in the leaf elements 
  in microseconds
\end{itemize}
If $\n{p}$ is not specified, only one process has been used.


\subsection{Robustness against mesh anisotropy, curvature and irregularity}
\label{sec:num-exp:robustness}

Previous studies have demonstrated the excellent convergence properties of 
the CS-MG-WAS method on uniform grids, but they have also highlighted its 
sensitivity to high aspect ratios
\cite{MG_Lottes2005a,MG_Stiller2016b,MG_Stiller2017a}.
Additional difficulties arise due to element deformations, irregular mesh structure
and variable coefficients.
These challenges are examined in more detail in the following subsections.
The test cases are based on Equation \eqref{eq:helmholtz:pde} with the exact 
solution
\begin{equation}
  \label{eq:num-exp:robustness:solution}
  u(\V{x}) = \sin(k x_1) \sin(k x_2) \sin(k x_3)
\end{equation}
and Dirichlet boundary conditions.
Unless stated otherwise, 
constant coefficients ${\kappa = 0}$, ${\nu = 1}$
and wave number ${k = 1}$
are assumed.
Depending on these parameters, the right-hand side $f(\V{x})$ is manufactured 
such that \eqref{eq:num-exp:robustness:solution} solves the problem with the given 
exact solution.

\begin{remark}
Since the tests focus on convergence in terms of residual reduction, the form of 
the exact solution has little influence. 
Studies using more complex solutions yield virtually identical results.
\end{remark}


\subsubsection{Uniform Cartesian}
\label{sec:num-exp:robustness:cart-uni}

For reference, a study was conducted using uniform Cartesian grids  
in the cubic domain ${[0,2\pi)^3}$ with periodic conditions.
The polynomial degree on the top level varied between $p=4$ and $32$.
To keep the problem size nearly constant, the number of elements was set to 
$128/p$ in each direction.
All multigrid methods introduced in Section~\ref{sec:mg-method} were tested
using the WAS and SCG smoothers with
${\n{s1} \le 16}$ steps for pre-smoothing and
${\n{s2} \le 8}$ for post-smoothing.
The FAS method used either CMG with ${\n{s0} = \n{s1} + \n{s2}}$ or 
FMG for starting and embedded interpolation or $L^2$ projection for 
fine-to-coarse solution projection.
$h$-multigrid was used with ${L=4}$ and 
$p$-multigrid with ${p_l = 2^{l-1}}$ up to ${p_L = p}$.
The multigrid methods were evaluated in more than 4000 tests using
penalties $\mu_{\star}$ from $2$ to $10$, 
overlaps $\delta$ from $0.04$ to $0.16$.
Each test was run on a single CPU core.
Table~\ref{tab:num-exp:robustness:cart-uni} 
compares the performance of different methods with 
fixed penalty ${\mu_{\star} = 2}$ and 
overlap ${\delta = 0.08}$.
Embedded interpolation was used for FAS solution projection.
First, we will focus on the $p$-MG methods shown on the left.
Before we begin comparing the methods, it should be noted that CS-MG-WAS and 
CS-CG-WAS were previously developed and analyzed specifically for the Cartesian case
\cite{MG_Stiller2017a}.
The current implementation achieves identical iteration counts on Cartesian grids,
but is also suitable for unstructured curvilinear meshes.
Among the $p$-MG methods, FAS-FMG-SCG is the consistently the fastest, although 
the optimal number of smoothing steps depends on the polynomial degree.
%
%
At the lowest degree, ${p=4}$, FAS proved to be less robust than CS when using 
the WAS smoother.
This drawback is eliminated by Krylov acceleration in SCG, which generally
improves the performance of $p$-MG.
With $h$-MG, Krylov acceleration improves the convergence rate only a slightly and 
therefore rarely reduces the number of iterations.
As a consequence, the cheaper WAS smoother achieves shorter runtimes and 
FAS-FMG-WAS becomes the fastest method.

Table~\ref{tab:num-exp:robustness:cart-uni:var} illustrates the influence of
various parameters on $p$-FAS-FMG-SCG and $h$-FAS-FMG-WAS.
In agreement with similar studies in \cite{MG_Fehn2020a}, we find that the methods 
are insensitive to the penalty parameter $\mu_{\star}$.
Using $L^2$ projection instead of embedded interpolation for $\NM P$ has no effect 
on $p$-FAS-FMG-SCG, but leads to a slight reduction of the convergence rate for
$h$-FAS-FMG-WAS.
Enlarging the overlap $\delta$ yields an increase in the convergence rate and, 
thus, can result in a reduction of the iteration count.
However, this advantage is offset by the higher cost of the Schwarz sweeps.
Finally, Table~\ref{tab:num-exp:robustness:cart-uni:fastest} summarizes 
the fastest runs for all tested parameters.
These results suggest that $h$-FAS-FMG-WAS is the most efficient method for 
solving problems on uniform Cartesian meshes, except for very high polynomial 
degrees where $p$-FAS-FMG-SCG is more suitable.
Although it is impossible to determine the optimum number of smoothing steps, 
${(\n{s1},\n{s2},\n{sc}) = (3,1,2)}$ is a good choice for $p$-FAS-FMG-SCG and
${(1,1,1)}$ for $h$-FAS-FMG-WAS.
It should also be noted that 
${\n{sc} = \lceil (\n{s1}+\n{s2})/2 \rceil}$
is generally a better choice for the number of crossover steps than
${\n{s1}+\n{s2}}$.

\begin{table}
\footnotesize
\begin{subtable}{0.48\textwidth}
\begin{NiceTabular}{lrcccrrc}
\toprule
 $p$-MG method & $p$ & $\n{s1}$ & $\n{s2}$ & $\n{sc}$ & $n_{10}$ 
 & $v_{10}$ & $\tau_{10}$ \\
\midrule

 CS-MG-WAS   &   4 &  1 &  1 &  2  & 10 & 10.0 & 5.10 \\
 CS-MG-WAS   &   4 &  2 &  2 &  4  &  6 & 10.0 & 5.00 \\
 CS-MG-SCG   &   4 &  2 &  2 &  4  &  4 & 12.0 & 4.78 \\
 CS-CG-WAS   &   4 &  2 &  2 &  4  &  5 &  9.9 & 4.05 \\
 FAS-CMG-WAS &   4 &  1 &  1 &  2  & 10 & 10.0 & 5.94 \\
 FAS-FMG-WAS &   4 &  1 &  1 &  2  & 10 & 10.0 & 6.07 \\
 FAS-FMG-WAS &   4 &  2 &  2 &  4  &  6 & 10.0 & 5.79 \\
 FAS-CMG-SCG &   4 &  2 &  2 &  4  &  3 &  9.0 & 4.07 \\
 FAS-FMG-SCG &   4 &  2 &  2 &  4  &  3 &  9.0 & 4.17 \\\RowStyle{\bfseries}
 FAS-FMG-SCG &   4 &  3 &  2 &  5  &  2 &  7.3 & 3.44 \\\midrule

 CS-MG-WAS   &   8 &  1 &  1 &  2  &  7 &  7.0 &  1.54 \\
 CS-MG-WAS   &   8 &  2 &  2 &  4  &  5 &  8.3 &  2.10 \\
 CS-MG-SCG   &   8 &  2 &  2 &  4  &  3 &  9.0 &  1.68 \\
 CS-CG-WAS   &   8 &  2 &  2 &  4  &  4 &  7.9 &  1.55 \\
 FAS-CMG-WAS &   8 &  1 &  1 &  2  &  9 &  9.0 &  2.48 \\
 FAS-FMG-WAS &   8 &  1 &  1 &  2  &  8 &  8.0 &  2.25 \\
 FAS-FMG-WAS &   8 &  2 &  2 &  4  &  6 & 10.0 &  2.62 \\
 FAS-CMG-SCG &   8 &  2 &  2 &  4  &  3 &  9.0 &  1.96 \\
 FAS-FMG-SCG &   8 &  2 &  2 &  4  &  3 &  9.0 &  2.04 \\\RowStyle{\bfseries}
 FAS-FMG-SCG &   8 &  3 &  2 &  3  &  2 &  6.1 &  1.52 \\\midrule

 CS-MG-WAS   &  16 &  1 &  1 &  2  &  6 &  6.0 &  0.94 \\
 CS-MG-WAS   &  16 &  2 &  2 &  4  &  4 &  6.7 &  1.03 \\
 CS-MG-SCG   &  16 &  2 &  2 &  4  &  3 &  9.0 &  1.32 \\
 CS-CG-WAS   &  16 &  2 &  2 &  4  &  3 &  5.9 &  0.82 \\
 FAS-CMG-WAS &  16 &  1 &  1 &  2  &  8 &  8.0 &  1.52 \\
 FAS-FMG-WAS &  16 &  1 &  1 &  2  &  7 &  7.0 &  1.36 \\
 FAS-FMG-WAS &  16 &  2 &  2 &  4  &  4 &  6.7 &  1.23 \\
 FAS-CMG-SCG &  16 &  2 &  2 &  4  &  3 &  9.0 &  1.37 \\
 FAS-FMG-SCG &  16 &  2 &  2 &  4  &  2 &  6.0 &  0.97 \\\RowStyle{\bfseries}
 FAS-FMG-SCG &  16 &  2 &  2 &  2  &  2 &  4.7 &  0.82 \\\midrule

 CS-MG-WAS   &  32 &  1 &  1 &  2  &  5 &  5.0 &  0.80 \\ 
 CS-MG-WAS   &  32 &  2 &  2 &  4  &  3 &  5.0 &  0.87 \\
 CS-MG-SCG   &  32 &  2 &  2 &  4  &  2 &  6.0 &  0.97 \\
 CS-CG-WAS   &  32 &  2 &  2 &  4  &  3 &  5.9 &  0.90 \\
 FAS-CMG-WAS &  32 &  1 &  1 &  2  &  8 &  8.0 &  1.64 \\
 FAS-FMG-WAS &  32 &  1 &  1 &  2  &  5 &  5.0 &  1.05 \\
 FAS-FMG-WAS &  32 &  2 &  2 &  4  &  3 &  5.0 &  1.02 \\
 FAS-CMG-SCG &  32 &  2 &  2 &  4  &  3 &  9.0 &  1.54 \\
 FAS-FMG-SCG &  32 &  2 &  2 &  4  &  2 &  6.0 &  1.09 \\\RowStyle{\bfseries}
 FAS-FMG-SCG &  32 &  3 &  1 &  2  &  1 &  1.7 &  0.55 \\

\bottomrule
\end{NiceTabular}

\caption{$p$-MG}
\label{tab:num-exp:robustness:cart-uni:p}
\end{subtable}
\hfill
\begin{subtable}{0.48\textwidth}
\begin{NiceTabular}{lrcccrrc}
\toprule
 $h$-MG method & $p$ & $\n{s1}$ & $\n{s2}$ & $\n{sc}$ & $n_{10}$ 
 & $v_{10}$ & $\tau_{10}$ \\
\midrule

 CS-MG-WAS   &   4 &  1 &   1 &  2 & 19 & 19.0 & 7.72 \\
 CS-MG-WAS   &   4 &  2 &   2 &  4 & 11 & 18.3 & 7.69 \\
 CS-MG-SCG   &   4 &  2 &   2 &  4 &  4 & 12.0 & 4.37 \\
 CS-CG-WAS   &   4 &  2 &   2 &  4 &  8 & 15.7 & 5.58 \\
 FAS-CMG-WAS &   4 &  1 &   1 &  2 &  9 &  9.0 & 4.19 \\
 FAS-FMG-WAS &   4 &  1 &   1 &  2 &  7 &  7.0 & 3.34 \\
 FAS-FMG-WAS &   4 &  2 &   2 &  4 &  4 &  6.7 & 3.26 \\
 FAS-CMG-SCG &   4 &  2 &   2 &  4 &  4 & 12.0 & 4.65 \\
 FAS-FMG-SCG &   4 &  2 &   2 &  4 &  3 &  9.0 & 3.54 \\\RowStyle{\bfseries}
 FAS-FMG-WAS &   4 &  3 &   1 &  2 &  3 &  3.8 & 1.93 \\\midrule

 CS-MG-WAS   &   8 &  1 &   1 &  2 &  8 &  8.0 & 1.67 \\
 CS-MG-WAS   &   8 &  2 &   2 &  4 &  6 & 10.0 & 2.02 \\
 CS-MG-SCG   &   8 &  2 &   2 &  4 &  4 & 12.0 & 2.26 \\
 CS-CG-WAS   &   8 &  2 &   2 &  4 &  5 &  9.8 & 1.79 \\
 FAS-CMG-WAS &   8 &  1 &   1 &  2 &  3 &  3.0 & 0.80 \\
 FAS-FMG-WAS &   8 &  1 &   1 &  2 &  1 &  1.0 & 0.35 \\
 FAS-FMG-WAS &   8 &  2 &   2 &  4 &  1 &  1.7 & 0.53 \\
 FAS-CMG-SCG &   8 &  2 &   2 &  4 &  2 &  6.0 & 1.30 \\
 FAS-FMG-SCG &   8 &  2 &   2 &  4 &  1 &  3.0 & 0.69 \\\RowStyle{\bfseries}
 FAS-FMG-WAS &   8 &  1 &   1 &  1 &  1 &  1.0 & 0.29 \\\midrule

 CS-MG-WAS   &  16 &  1 &   1 &  2 &  9 &  9.0 & 1.62 \\
 CS-MG-WAS   &  16 &  2 &   2 &  4 &  6 & 10.0 & 1.72 \\
 CS-MG-SCG   &  16 &  2 &   2 &  4 &  4 & 12.0 & 1.97 \\
 CS-CG-WAS   &  16 &  2 &   2 &  4 &  5 &  9.8 & 1.57 \\
 FAS-CMG-WAS &  16 &  1 &   1 &  2 &  3 &  3.0 & 0.71 \\
 FAS-FMG-WAS &  16 &  1 &   1 &  2 &  2 &  2.0 & 0.53 \\
 FAS-FMG-WAS &  16 &  2 &   2 &  4 &  2 &  3.3 & 0.76 \\
 FAS-CMG-SCG &  16 &  2 &   2 &  4 &  2 &  6.0 & 1.13 \\
 FAS-FMG-SCG &  16 &  2 &   2 &  4 &  2 &  6.0 & 1.14 \\\RowStyle{\bfseries}
 FAS-FMG-WAS &  16 &  1 &   1 &  1 &  2 &  1.7 & 0.44 \\\midrule

 CS-MG-WAS   &  32 &  1 &   1 &  2 &  4 &  4.0 & 1.18 \\ 
 CS-MG-WAS   &  32 &  2 &   2 &  4 &  3 &  5.0 & 1.30 \\
 CS-MG-SCG   &  32 &  2 &   2 &  4 &  3 &  9.0 & 1.89 \\
 CS-CG-WAS   &  32 &  2 &   2 &  4 &  3 &  5.9 & 1.31 \\
 FAS-CMG-WAS &  32 &  1 &   1 &  2 &  3 &  3.0 & 1.11 \\
 FAS-FMG-WAS &  32 &  1 &   1 &  2 &  3 &  3.0 & 1.13 \\
 FAS-FMG-WAS &  32 &  2 &   2 &  4 &  2 &  3.3 & 1.08 \\
 FAS-CMG-SCG &  32 &  2 &   2 &  4 &  2 &  6.0 & 1.45 \\
 FAS-FMG-SCG &  32 &  2 &   2 &  4 &  2 &  6.0 & 1.44 \\\RowStyle{\bfseries}
 FAS-FMG-WAS &  32 &  2 &   2 &  2 &  2 &  2.7 & 0.89 \\

\bottomrule
\end{NiceTabular}
\caption{$h$-MG}
\label{tab:num-exp:robustness:cart-uni:h}
\end{subtable}
\caption{Uniform Cartesian mesh: performance of selected methods with
${\mu_{\star} = 2}$ and ${\delta = 0.08}$.
FAS methods used embedded interpolation for projection
and FAS-CMG ${\n{s0}=\n{s1}+\n{s2}}$ iterations in the starting cascade.
Fastest runs are shown in bold.}
\label{tab:num-exp:robustness:cart-uni}
\end{table}

\begin{table}
\footnotesize
\begin{subtable}{0.48\textwidth}
\begin{NiceTabular}{rrcccccrrc}
\toprule
 $p$ & $\mu_\star$ & $\delta$ & $\protect\underaccent{\bar}{P}$ 
     & $\n{s1}$ & $\n{s2}$ & $\n{sc}$ & $n_{10}$ 
     & $v_{10}$ & $\tau_{10}$ \\
\midrule

  4  &   2  &  0.08  &  I  &   3  &  2  &  5  &  2 & 7.3 & 3.44 \\
  4  &   2  &  0.08  &  L  &   3  &  2  &  5  &  2 & 7.3 & 3.44 \\
  4  &   5  &  0.08  &  I  &   3  &  2  &  5  &  2 & 7.3 & 3.45 \\
  4  &  10  &  0.08  &  I  &   3  &  2  &  5  &  2 & 7.3 & 3.44 \\\midrule

  8  &   2  &  0.08  &  I  &   3  &  2  &  3  &  2 & 6.1 & 1.52 \\
  8  &   2  &  0.08  &  L  &   3  &  2  &  3  &  2 & 6.1 & 1.46 \\
  8  &   5  &  0.08  &  I  &   3  &  2  &  3  &  2 & 6.1 & 1.52 \\
  8  &  10  &  0.08  &  I  &   3  &  2  &  3  &  2 & 6.1 & 1.53 \\
  8  &   2  &  0.04  &  I  &   3  &  2  &  3  &  3 & 9.7 & 2.07 \\\midrule

 16  &   2  &  0.08  &  I  &   2  &  2  &  2  &  2 & 4.7 & 0.82 \\
 16  &   2  &  0.08  &  L  &   2  &  2  &  2  &  2 & 4.7 & 0.83 \\
 16  &   5  &  0.08  &  I  &   2  &  2  &  2  &  2 & 4.7 & 0.92 \\
 16  &  10  &  0.08  &  I  &   2  &  2  &  2  &  2 & 4.7 & 0.87 \\
 16  &   2  &  0.04  &  I  &   2  &  2  &  2  &  3 & 7.7 & 1.17 \\
 16  &   2  &  0.16  &  I  &   2  &  2  &  2  &  2 & 4.7 & 1.05 \\\midrule

 32  &   2  &  0.08  &  I  &   3  &  1  &  2  &  1 & 1.7 & 0.55 \\
 32  &   2  &  0.08  &  L  &   3  &  1  &  2  &  1 & 1.7 & 0.55 \\
 32  &   5  &  0.08  &  I  &   3  &  1  &  2  &  1 & 1.7 & 0.56 \\
 32  &  10  &  0.08  &  I  &   3  &  1  &  2  &  1 & 1.7 & 0.57 \\
 32  &   2  &  0.04  &  I  &   3  &  1  &  2  &  1 & 1.7 & 0.50 \\
 32  &   2  &  0.16  &  I  &   3  &  1  &  2  &  1 & 1.7 & 0.70 \\

\bottomrule
\end{NiceTabular}

\caption{$p$-FAS-FMG-SCG}
\label{tab:num-exp:robustness:cart-uni:p:var}
\end{subtable}
\hfill
\begin{subtable}{0.48\textwidth}
\begin{NiceTabular}{rrcccccrrc}
\toprule
 $p$ & $\mu_\star$ & $\delta$ & $\protect\underaccent{\bar}{P}$ 
     & $\n{s1}$ & $\n{s2}$ & $\n{sc}$ & $n_{10}$ 
     & $v_{10}$ & $\tau_{10}$ \\
\midrule

  4  &   2  &  0.08  &  I  &  3  &  1  &  2  &  3  &  3.8 & 1.93 \\
  4  &   2  &  0.08  &  L  &  3  &  1  &  2  &  3  &  3.8 & 1.93 \\
  4  &   5  &  0.08  &  I  &  3  &  1  &  2  &  3  &  3.8 & 1.93 \\
  4  &  10  &  0.08  &  I  &  3  &  1  &  2  &  3  &  3.8 & 1.93 \\\midrule
     
  8  &   2  &  0.08  &  I  &  1  &  1  &  1  &  1  &  1.0 & 0.29 \\
  8  &   2  &  0.08  &  L  &  1  &  1  &  1  &  1  &  1.0 & 0.29 \\
  8  &   5  &  0.08  &  I  &  1  &  1  &  1  &  1  &  1.0 & 0.29 \\
  8  &  10  &  0.08  &  I  &  1  &  1  &  1  &  1  &  1.0 & 0.29 \\
  8  &   2  &  0.04  &  I  &  1  &  1  &  1  &  5  &  3.8 & 0.83 \\\midrule
     
 16  &   2  &  0.08  &  I  &  1  &  1  &  1  &  2  &  1.7 & 0.44 \\
 16  &   2  &  0.08  &  L  &  1  &  1  &  1  &  3  &  2.4 & 0.58 \\
 16  &   5  &  0.08  &  I  &  1  &  1  &  1  &  3  &  2.4 & 0.60\\
 16  &  10  &  0.08  &  I  &  1  &  1  &  1  &  3  &  2.4 & 0.59 \\
 16  &   2  &  0.04  &  I  &  1  &  1  &  1  &  4  &  3.1 & 0.68 \\
 16  &   2  &  0.16  &  I  &  1  &  1  &  1  &  1  &  1.0 & 0.35 \\\midrule
     
 32  &   2  &  0.08  &  I  &  2  &  2  &  2  &  3  &  2.4 & 0.97 \\
 32  &   2  &  0.08  &  L  &  2  &  2  &  2  &  3  &  2.4 & 0.94 \\
 32  &   5  &  0.08  &  I  &  2  &  2  &  2  &  3  &  2.4 & 0.96 \\
 32  &  10  &  0.08  &  I  &  2  &  2  &  2  &  2  &  2.7 & 0.86 \\
 32  &   2  &  0.04  &  I  &  2  &  2  &  2  &  4  &  3.1 & 1.28 \\
 32  &   2  &  0.16  &  I  &  2  &  2  &  2  &  2  &  2.7 & 1.01 \\

\bottomrule
\end{NiceTabular}

\caption{$h$-FAS-FMG-WAS}
\label{tab:num-exp:robustness:cart-uni:h:var}
\end{subtable}
\caption{Uniform Cartesian mesh: FAS-FMG performance with
varying penalties $\mu_{\star}$ and overlaps $\delta$
including a different number of node layers.
${\protect\underaccent{\bar}{P} = \text{I}}$
denotes solution projection with embedded interpolation, and
L with $L^2$-projection.
}
\label{tab:num-exp:robustness:cart-uni:var}
\end{table}

\begin{table}
\footnotesize
\begin{subtable}{0.49\textwidth}
\begin{NiceTabular}{rrcccccrrc}
\toprule
 $p$ & $\mu_\star$ & $\delta$ & $\protect\underaccent{\bar}{P}$ 
     & $\n{s1}$ & $\n{s2}$ & $\n{sc}$ & $n_{10}$ 
     & $v_{10}$ & $\tau_{10}$ \\
\midrule

  4  &   5  &  0.08  &  I  &   3  &  1  &  2  &  2  &  5.5  & 2.66 \\
  8  &  10  &  0.04  &  I  &   3  &  1  &  2  &  2  &  5.5  & 1.23 \\
 16  &  10  &  0.04  &  I  &   3  &  1  &  2  &  1  &  1.8  & 0.60 \\
 32  &   5  &  0.04  &  I  &   2  &  2  &  2  &  1  &  1.7  & 0.43 \\

\bottomrule
\end{NiceTabular}

\caption{$p$-FAS-FMG-SCG}
\label{tab:num-exp:robustness:cart-uni:p:fastest}
\end{subtable}
\hfill
\begin{subtable}{0.49\textwidth}
\begin{NiceTabular}{rrcccccrrc}
\toprule
 $p$ & $\mu_\star$ & $\delta$ & $\protect\underaccent{\bar}{P}$ 
     & $\n{s1}$ & $\n{s2}$ & $\n{sc}$ & $n_{10}$ 
     & $v_{10}$ & $\tau_{10}$ \\
\midrule

  4  &  10  &  0.08  &  I  &  3  &  2  &  3  &  1  &  3.4  & 1.76 \\
  8  &   5  &  0.08  &  I  &  1  &  1  &  1  &  1  &  1.0  & 0.29 \\
 16  &   5  &  0.16  &  I  &  1  &  1  &  1  &  1  &  1.0  & 0.34 \\
 32  &   5  &  0.16  &  I  &  1  &  1  &  1  &  2  &  1.7  & 0.78 \\

\bottomrule
\end{NiceTabular}

\caption{$h$-FAS-FMG-WAS}
\label{tab:num-exp:robustness:cart-uni:h:fastest}
\end{subtable}
\caption{Uniform Cartesian mesh: fastest methods.
For caption see Table~\ref{tab:num-exp:robustness:cart-uni:var}.}
\label{tab:num-exp:robustness:cart-uni:fastest}
\end{table}


\subsubsection{Cartesian with varying aspect ratio}
\label{sec:num-exp:robustness:cart-var}

\begin{figure}
  \subcaptionbox{$p$-MG iteration counts
    \label{fig:num-exp:cart-var:p:n10}}
    {\includegraphics[scale=0.52]{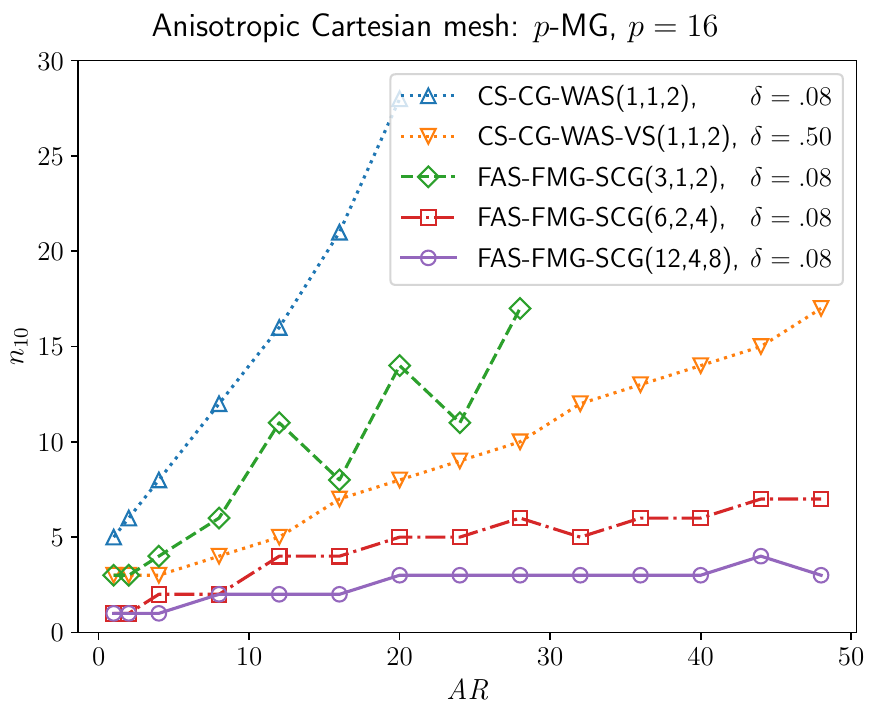}}
  \hfill
  \subcaptionbox{$h$-MG iteration counts
    \label{fig:num-exp:cart-var:h:n10}}
    {\includegraphics[scale=0.52]{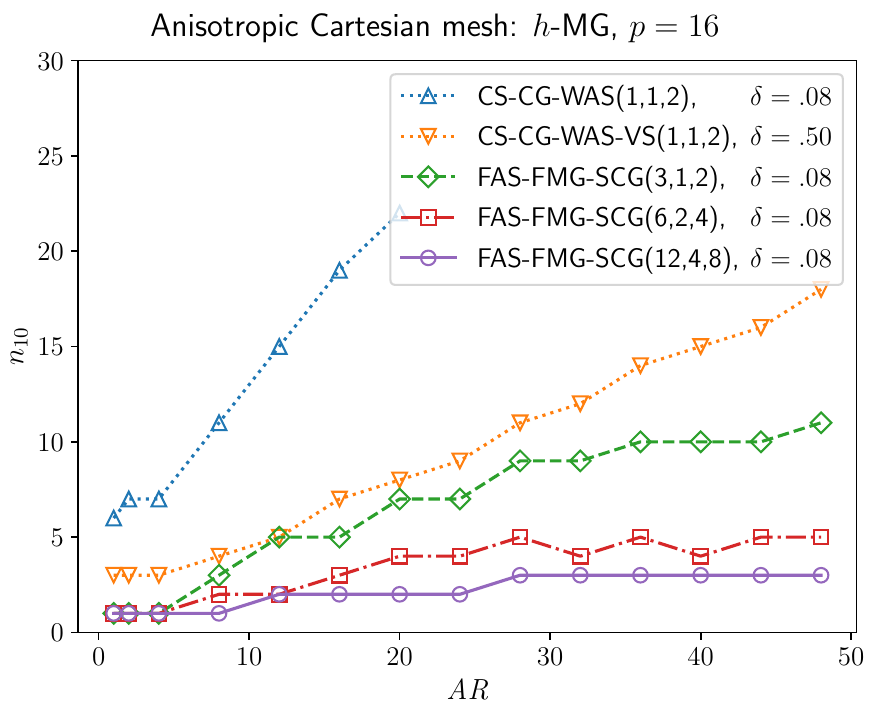}}
  \\[\medskipamount]
  \subcaptionbox{$p$-MG normalized runtime
    \label{fig:num-exp:cart-var:p:tau10}}
    {\includegraphics[scale=0.52]{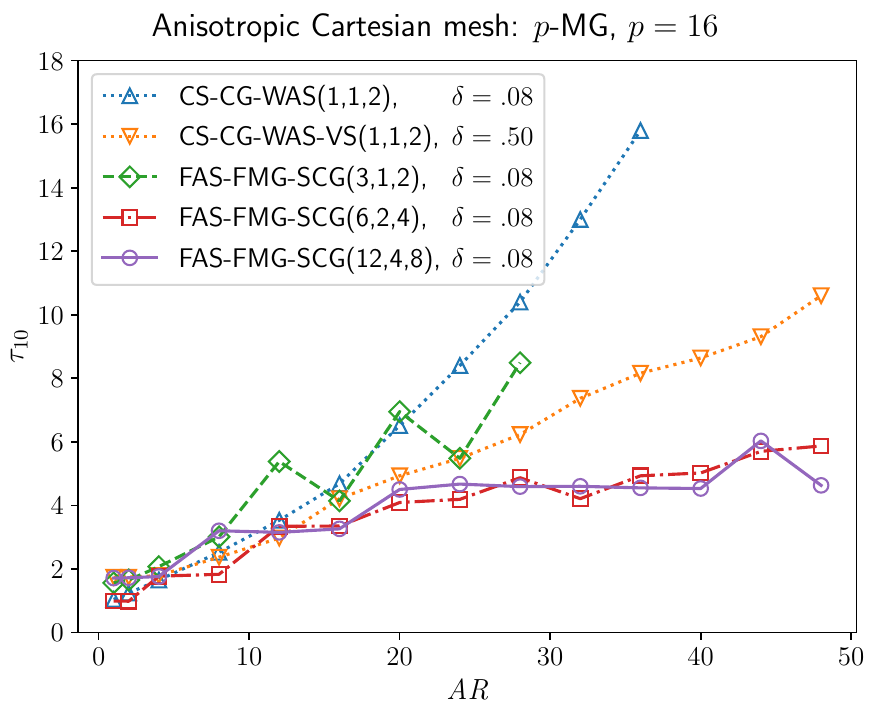}}
  \hfill
  \subcaptionbox{$h$-MG normalized runtime
    \label{fig:num-exp:cart-var:h:tau10}}
    {\includegraphics[scale=0.52]{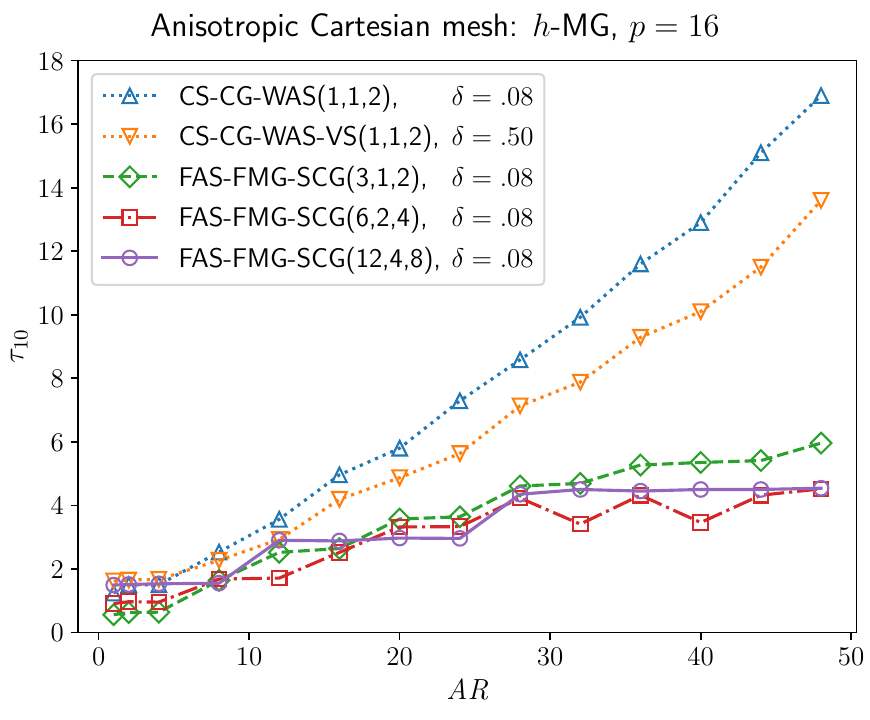}}
  \caption{Robustness of MG methods with respect to increasing element
    aspect ratios
    \label{fig:num-exp:cart-var}}
\end{figure}

The sensitivity of MG methods to the element aspect ratio is a well-known issue.
As a countermeasure, \citet{MG_Lottes2005a} enhanced the robustness by using 
Krylov acceleration on the top level.
This approach was further improved in
\cite{MG_Stiller2017a}
by employing larger subdomain overlaps and variable smoothing.
\citet{MG_Fehn2020a} applied combinations of $p$- and $h$-MG methods
with Chebyshev smoothers to test cases with different aspect ratios.
Although the authors have not presented a detailed study on this issue,
their findings suggest that the iteration count increases by a factor of 
5 to 7 as the aspect ratio approaches 67.
\citet{MG_Phillips2025a} achieved comparable results with $p$-MG methods
using optimized Jacobi smoothers and Krylov acceleration on the top level.
To examine the robustness of the methods developed in the present study, we modify 
the test case presented in Section~\ref{sec:num-exp:robustness:cart-uni} by stretching 
the computational domain in the first and second directions such that 
${\Omega = [0, 2\pi \mathit{AR}) \times
           [0, 2\pi \lceil\mathit{AR}/2\rceil) \times
           [0, 2\pi)}$.
The multilevel meshes are constructed as in the uniform case, such that 
$\mathit{AR}$ represents the ratio of the longest to the shortest element
side.
Figure~\ref{fig:num-exp:cart-var} presents the iteration counts $n_{10}$ and 
runtimes $\tau_{10}$ of selected MG methods for ${p = 16}$ in the range
${1 \le \mathit{AR} \le 48}$.
Methods using the WAS smoother without Krylov acceleration are not shown,
as they become soon inefficient or even unstable as the aspect ratio increases.
Initially, we observe that the accelerated correction scheme CS-CG-WAS(1,1,2), 
with an 8 percent overlap, succeeds but exhibits a rapidly escalating iteration 
count.
This results in prohibitive computational costs for AR approaching 48.
In accordance with earlier studies, the number of iterations can be reduced
significantly by enlarging the overlap and variable smoothing.
This is illustrated using the example of CS-CG-WAS-VS(1,1,2), 
which employs an overlap of ${\delta = 0.5}$ and
twice the number of smoothing steps on each coarser level, i.e.
${\n[,l]{s1} = 2 \n[,l+1]{s1}}$ etc.
Unfortunately, these measures entail additional computational costs,
which diminish the gain in runtime.
The MG methods using the SCG smoother benefit from the Krylov acceleration 
applied on each level and are therefore much more robust than their MG-WAS 
counterparts.
Unfortunately, the small number of smoothing steps leading to the fastest 
methods on uniform grids proved insufficient to maintain robustness 
against large aspect ratios.
However, increasing the number of steps yields extremely robust 
and still very fast methods.
For example, FAS-FMG-SCG(12,4,8) achieves a nearly constant iteration count
${n_{10} \le 4}$ and normalized runtimes ${\tau_{10} < 5}$ over the whole range 
of aspect ratios.
As shown in Table~\ref{tab:num-exp:robustness:cart-var:fastest},
even faster methods are available, depending on the polynomial 
degree, the aspect ratio and the coarsening method.

\begin{table}
\footnotesize
\begin{subtable}{0.45\textwidth}
\begin{NiceTabular}{rrcccccrrc}
\toprule
 $p$ & $\mathit{AR}$ & $\delta$  
     & $\n{s1}$ & $\n{s2}$ & $\n{sc}$ & $n_{10}$ 
     & $v_{10}$ & $\tau_{10}$ \\
\midrule

  4  &  12  &  0.24  &  12  &  4  &  8  &  3  &  28.0  & 18.6 \\
  4  &  48  &  0.04  &  12  &  4  &  8  &  8  &  83.0  & 71.7 \\
\midrule
  8  &  12  &  0.24  &   8  &  2  &  5  &  3  &  17.9  & 5.51 \\
  8  &  48  &  0.24  &  12  &  4  &  8  &  3  &  28.0  & 9.16 \\
\midrule
 16  &  12  &  0.24  &  12  &  4  &  8  &  1  &   6.0  & 2.46 \\
 16  &  48  &  0.24  &  12  &  4  &  8  &  2  &  17.0  & 4.57 \\
\midrule
 32  &  12  &  0.24  &   6  &  2  &  4  &  1  &   3.2  & 1.50 \\
 32  &  48  &  0.08  &  12  &  4  &  8  &  2  &  17.0  & 3.23 \\

\bottomrule
\end{NiceTabular}

\caption{$p$-FAS-FMG-SCG}
\label{tab:num-exp:robustness:cart-var:p:fastest}
\end{subtable}
\qquad
\begin{subtable}{0.45\textwidth}
\begin{NiceTabular}{rrcccccrrc}
\toprule
 $p$ & $\mathit{AR}$ & $\delta$  
     & $\n{s1}$ & $\n{s2}$ & $\n{sc}$ & $n_{10}$ 
     & $v_{10}$ & $\tau_{10}$ \\
\midrule

  4  &  12  &  0.24  &   8  &  2  &  5  &  3  &  17.9  & 6.73 \\
  4  &  48  &  0.24  &   6  &  2  &  4  &  7  &  37.2  & 12.6 \\
\midrule
  8  &  12  &  0.24  &   8  &  2  &  5  &  2  &  10.9  & 2.95 \\
  8  &  48  &  0.08  &   8  &  2  &  5  &  5  &  31.9  & 5.59 \\
\midrule
 16  &  12  &  0.24  &   3  &  1  &  2  &  2  &   4.7  & 1.52 \\
 16  &  48  &  0.24  &   8  &  2  &  5  &  2  &  10.9  & 3.21 \\
\midrule
 32  &  12  &  0.08  &   3  &  1  &  2  &  2  &   4.7  & 2.56 \\
 32  &  48  &  0.24  &   6  &  2  &  4  &  2  &   8.8  & 6.07 \\

\bottomrule
\end{NiceTabular}

\caption{$h$-FAS-FMG-SCG}
\label{tab:num-exp:robustness:cart-var:h:fastest}
\end{subtable}
\caption{Cartesian mesh with varying element aspect ratio $\mathit{AR}$: 
fastest methods using ${\mu_\star=5}$ and embedded interpolation for projection}
\label{tab:num-exp:robustness:cart-var:fastest}
\end{table}


\subsubsection{Curved}
\label{sec:num-exp:robustness:curved}

\begin{figure}
  \includegraphics[height=0.4\textheight]{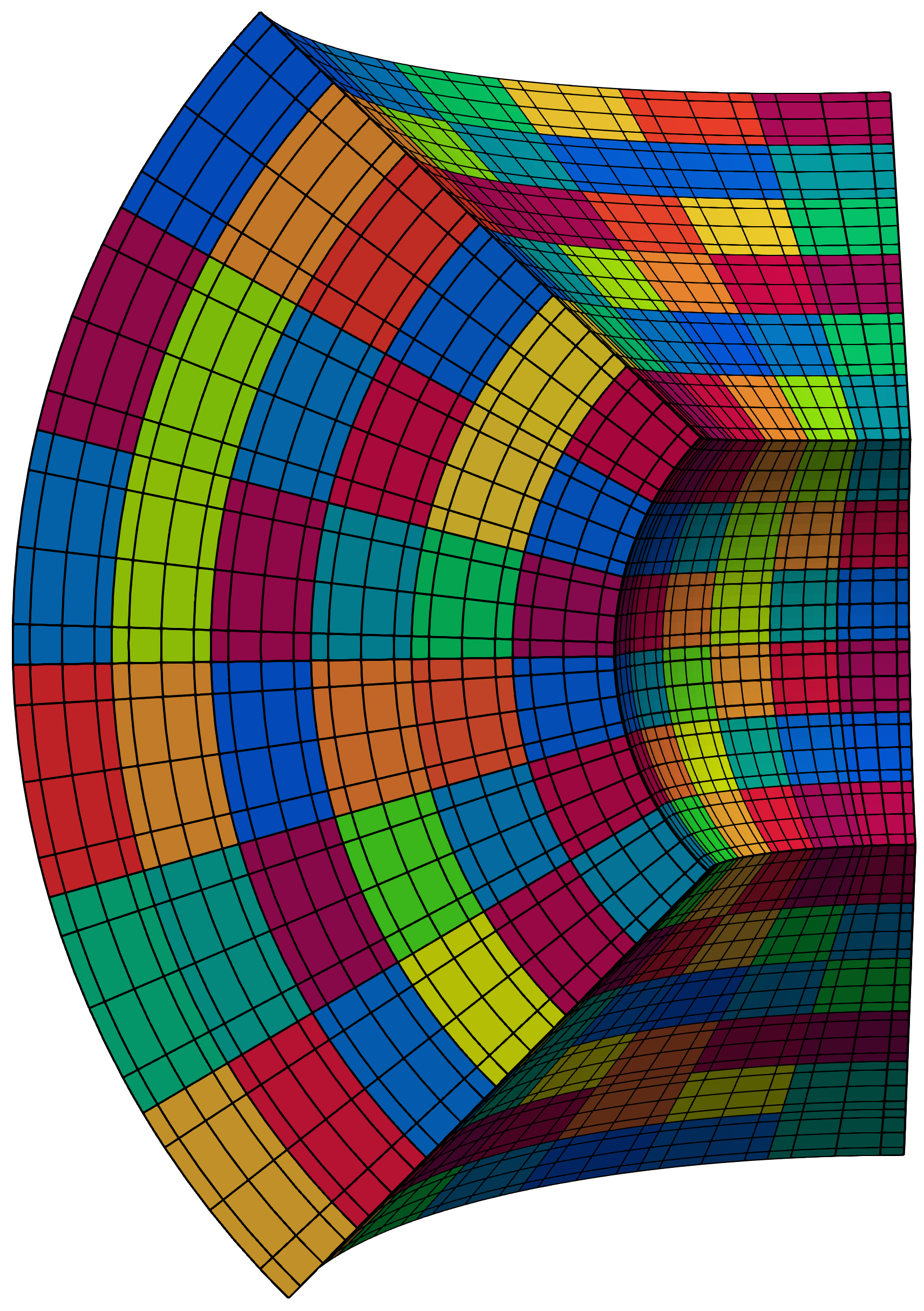}
  \caption{Spherical shell segment with $6^3$ elements of degree 4.
    Colors identify individual elements and lines the mesh spanned by 
    collocation points.
    \label{fig:num-exp:curved:mesh}}
\end{figure}

To examine the  robustness against curvature, the problem is considered 
in an octant of a spherical shell.
The computational domain is given in spherical coordinates ${(r,\theta,\varphi)}$ by
\begin{equation*}
  \Omega = \left( \frac{1}{2}  , \frac{3}{2}   \right) \times
           \left(-\frac{\pi}{4}, \frac{\pi}{4} \right) \times
           \left(-\frac{\pi}{4}, \frac{\pi}{4} \right)
           \,.
\end{equation*}
For illustration, Figure~\ref{fig:num-exp:curved:mesh} shows a decomposition of the domain
into $6^3$ elements of degree 4.
Numerical experiments were conducted using $(128/p)^3$ elements 
with ${p \in \{4,8,16,32\}}$.
This corresponds to element lengths of 
${h_r = p / 128}$
in the radial direction and
${h_{\theta} = h_{\varphi} = r \pi p / 256}$ 
in the meridional and azimuthal directions.
Consequently, the normalized curvature 
${\kappa_{\mathrm c} = h_{\varphi} / r}$
varies from ${\pi/64}$ with ${p=4}$ to ${\pi/8}$ with ${p=32}$. 
The aspect ratio is limited to ${3\pi/4}$, regardless of the degree $p$.
Table~\ref{tab:num-exp:robustness:curved:fastest} shows the results of
the fastest $p$-FAS-FMG-SCG and $h$-FAS-FMG-WAS methods within 
the tested parameter range.
Compared to uniform Cartesian grids
(Table~\ref{tab:num-exp:robustness:cart-uni:fastest}), 
more smoothing steps are applied within one V-cycle.
This results in a moderate increase in runtime, 
typically by about a factor of two.
Only in the case of $p$-FAS-FMG-SCG with ${p=4}$ is the factor of three exceeded.

\begin{table}
\footnotesize
\begin{subtable}{0.48\textwidth}
\begin{NiceTabular}{rrcccccrrc}
\toprule
 $p$ & $\mu_\star$ & $\delta$ & $\protect\underaccent{\bar}{P}$ 
     & $\n{s1}$ & $\n{s2}$ & $\n{sc}$ & $n_{10}$ 
     & $v_{10}$ & $\tau_{10}$ \\
\midrule

  4  &  5  &  0.08  &  I  &  8  &  2  &  5  &  2  &  10.9  &  9.28  \\ 
  8  &  5  &  0.04  &  I  &  6  &  2  &  4  &  2  &   8.8  &  2.80  \\ 
 16  &  5  &  0.04  &  I  &  6  &  2  &  4  &  1  &   5.7  &  1.00  \\ 
 32  &  5  &  0.04  &  I  &  3  &  1  &  2  &  1  &   3.0  &  0.71  \\ 

\bottomrule
\end{NiceTabular}

\caption{$p$-FAS-FMG-SCG}
\label{tab:num-exp:robustness:curved:p:fastest}
\end{subtable}
\hfill
\begin{subtable}{0.48\textwidth}
\begin{NiceTabular}{rrcccccrrc}
\toprule
 $p$ & $\mu_\star$ & $\delta$ & $\protect\underaccent{\bar}{P}$ 
     & $\n{s1}$ & $\n{s2}$ & $\n{sc}$ & $n_{10}$ 
     & $v_{10}$ & $\tau_{10}$ \\
\midrule


   4  &  5  &  0.08  &  I  &  6  &  2  &  4  &  1  &  3.0  & 2.14 \\ %
   8  &  5  &  0.04  &  I  &  3  &  1  &  2  &  1  &  1.7  & 0.52 \\ 
  16  &  5  &  0.08  &  I  &  3  &  1  &  2  &  1  &  1.7  & 0.51 \\ 
  32  &  5  &  0.04  &  I  &  6  &  2  &  4  &  1  &  3.0  & 1.12 \\ 

\bottomrule
\end{NiceTabular}

\caption{$h$-FAS-FMG-SCG}
\label{tab:num-exp:robustness:curved:h:fastest}
\end{subtable}
\caption{Performance of selected $p$-MG and $h$-MG methods on the curved mesh
  test case.}
\label{tab:num-exp:robustness:curved:fastest}
\end{table}


\subsubsection{Irregular}
\label{sec:num-exp:robustness:diamond}

\begin{figure}
  \subcaptionbox{2D unit cell
    \label{fig:num-exp:robustness:diamond:mesh:cell}}
    {\hspace*{1em}
     \includegraphics[height=0.15\textheight]{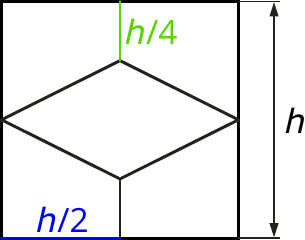}
     \vspace*{1\baselineskip}}
  \qquad\quad
  \subcaptionbox{
    6 layers of 3$\times$3 unit cells with elements of degree 4
    \label{fig:num-exp:uns-diamond:mesh:example}}
    {\includegraphics[height=0.4\textheight]{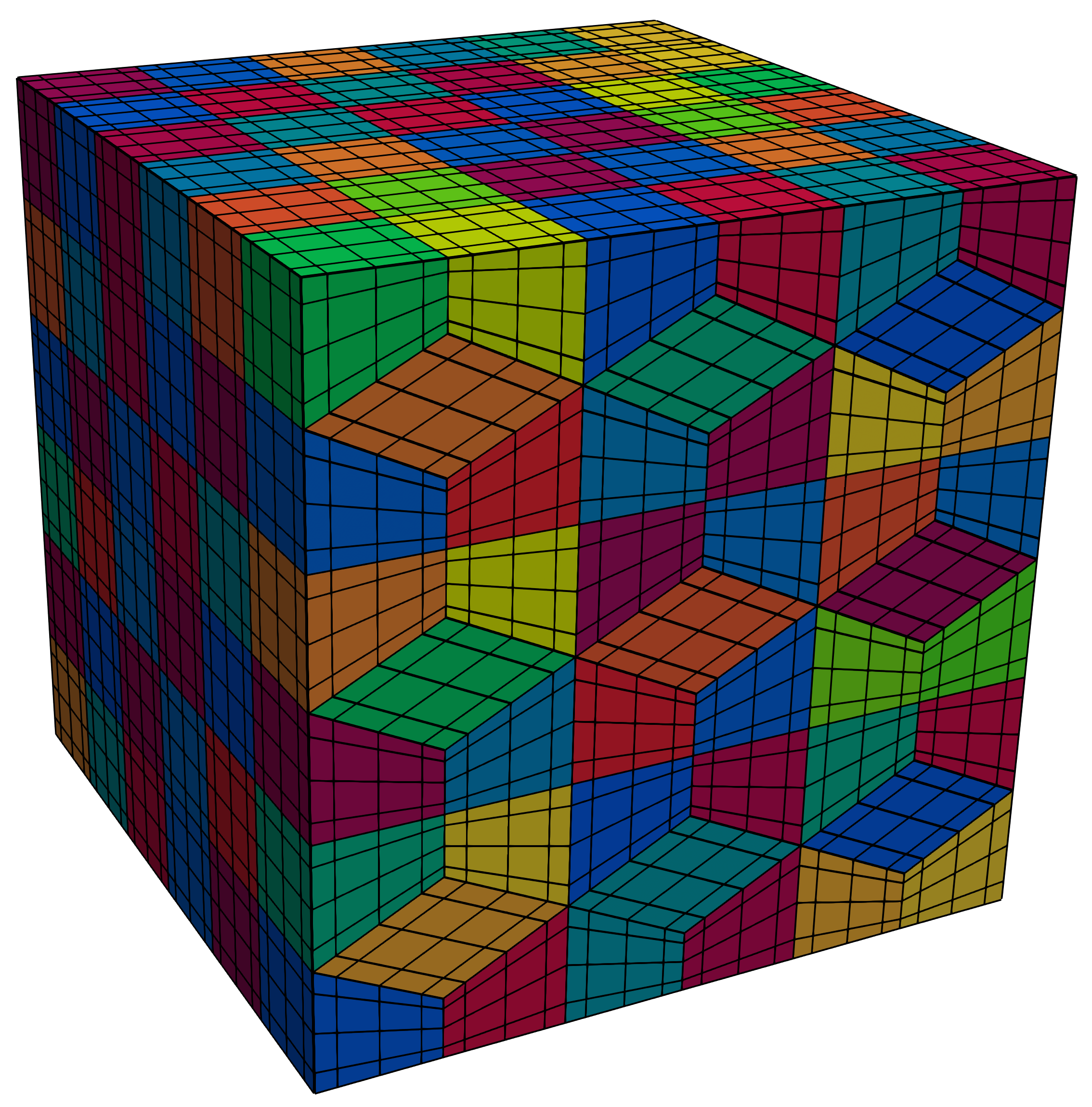}}
  \caption{Diamond mesh: on the left, a square unit cell with a rhomboid inset; 
     on the right, a mesh with elements of degree 4, comprising 6 layers of  
     3$\times$3 cells
    \label{fig:num-exp:robustness:diamond:mesh}}
\end{figure}

Next, we revisit the problem examined in 
Section~\ref{sec:num-exp:robustness:cart-uni}
on an irregular mesh.
The mesh is constructed from 2D unit cells, which are replicated in the
$x_1$-$x_2$-plane and then extruded in the $x_3$-direction to generate 
the requested layers of elements
(Figure~\ref{fig:num-exp:robustness:diamond:mesh}),
Each unit cell forms a square that is divided into four trapezoids at its 
corners and into a rhombus (diamond) in the center.
This structure yields a mesh in which about one third of the edges and 
two thirds of the vertices are are shared by either more or fewer elements 
than in the Cartesian case.
The resulting subdomains are severely distorted and irregular, 
which renders the tensor-product Schwarz method inefficient.
Consequently, the MG methods based on the WAS smoother fail to converge.
In contrast, $p$-FAS-FMG-SCG is only moderately affected %
(Table~\ref{tab:num-exp:robustness:diamond:p:fastest}).
The $h$-FAS-FMG-SCG method proves to be even more robust, suffering only 
a minor performance loss compared to the uniform Cartesian mesh
(Table~\ref{tab:num-exp:robustness:diamond:h:fastest}).
Doubling the overlap to ${\delta=0.08}$ results in the same number of 
iterations, but leads to a slightly longer runtime.
The inclusion of additional neighbor nodes at irregular edges or corners, 
as suggested in Section ~\ref{sec:mg-method:smoothers:schwarz}, 
resulted in no significant improvement, except for a few marginally stable 
cases using a very small number of smoothing steps.

\begin{table}
\footnotesize
\begin{subtable}{0.48\textwidth}
\begin{NiceTabular}{rrcccccrrc}
\toprule
 $p$ & $\mu_\star$ & $\delta$ & $\protect\underaccent{\bar}{P}$ 
     & $\n{s1}$ & $\n{s2}$ & $\n{sc}$ & $n_{10}$ 
     & $v_{10}$ & $\tau_{10}$ \\
\midrule

 
   4  &  5  &  0.04  &  I  &   3  &  1  &  2  &  4  &   8.2  &  12.1  \\
   8  &  5  &  0.04  &  I  &   8  &  2  &  5  &  2  &  10.9  &  3.63  \\
  16  &  5  &  0.04  &  I  &   8  &  2  &  5  &  1  &   7.0  &  1.21  \\
  32  &  5  &  0.04  &  I  &   8  &  2  &  5  &  1  &   7.0  &  1.15  \\
  
  

\bottomrule
\end{NiceTabular}

\caption{$p$-FAS-FMG-SCG}
\label{tab:num-exp:robustness:diamond:p:fastest}
\end{subtable}
\hfill
\begin{subtable}{0.48\textwidth}
\begin{NiceTabular}{rrcccccrrc}
\toprule
 $p$ & $\mu_\star$ & $\delta$ & $\protect\underaccent{\bar}{P}$ 
     & $\n{s1}$ & $\n{s2}$ & $\n{sc}$ & $n_{10}$ 
     & $v_{10}$ & $\tau_{10}$ \\
\midrule

 
 
 
   4  &  5  &  0.04  &  I  &   3  &  1  &  2  &  3  &   6.7  &  4.02  \\
   8  &  5  &  0.04  &  I  &   6  &  2  &  4  &  1  &   5.7  &  1.06  \\
  16  &  5  &  0.04  &  I  &   3  &  1  &  2  &  1  &   3.0  &  0.62  \\
  32  &  5  &  0.04  &  I  &   3  &  1  &  2  &  1  &   3.0  &  0.81  \\

\bottomrule
\end{NiceTabular}

\caption{$h$-FAS-FMG-SCG}
\label{tab:num-exp:robustness:diamond:h:fastest}
\end{subtable}
\caption{Performance of selected $p$-MG and $h$-MG methods in diamond mesh 
  test case.}
\label{tab:num-exp:robustness:diamond:fastest}
\end{table}


\subsubsection{Blade}
\label{sec:num-exp:robustness:blade}

The robustness study concludes with a test case derived from the analysis of 
flow in turbo compressors.
Here we consider a simplified configuration representing a linear cascade 
with a single blade.
Figure~\ref{fig:num-exp:robustness:blade:mesh:coarse} shows the computational
domain and the initial single-level mesh comprising 20640 elements.
This mesh combines all the challenges examined before: 
aspect ratios of up to 10, 
curved and distorted elements
as well as
irregular edges and vertices.
Additionally,
the element size varies considerable
due to the refinement toward the blade.
To mimic the implicit viscous step of a semi-implicit flow-solver, 
Equation~\ref{eq:helmholtz:pde} is considered
with 
${\kappa = 1}$,
solution
${u = \sin(k x_1) \sin(k x_2) \sin(k x_3)}$,
variable diffusivity
${\nu = 1 + \cos(k x_1) \cos(k x_2) \cos(k x_3)}$,
wave number
${k = 2\pi}$
and 
Dirichlet conditions at all boundaries.
For examining the performance of the FAS-FMG-SCG method, we define 
one test case with global refinement and another with local refinement.
In both cases, the root level is identical to the initial mesh with ${p = 1}$.
The multilevel mesh is constructed by applying a sequence of 
three $p$-refinements and two $h$-refinements.
In the global case, the top level contains 1320960 elements of degree 8, 
resulting in 963 million degrees of freedom.
In the local case, only elements adjacent to the blade are refined.
This yields a multilevel mesh comprising 79440 leaf elements with
46 million degrees of freedom, which are distributed over five levels.
For illustration, 
Figure~\ref{fig:num-exp:robustness:blade:mesh:refined} 
shows the composite mesh near the blade.
Table~\ref{tab:num-exp:robustness:blade:mesh-lr}
summarizes the main mesh characteristics.
For assessing mesh quality, we consider the maximum aspect ratio 
\begin{equation*}
  \mathit{AR}_{\max}
  = \max_e \frac{\Delta x^e_{\max}}{\Delta x^e_{\min}}
\end{equation*}
based on the minimum and maximum extensions of the surrogate element domains,
(see Section~\ref{sec:mg-method:smoothers:schwarz})
and
the minimum of the scaled element Jacobians
\begin{equation*}
  Q_{\min} = \min_e \frac{J^e_{\min}}{J^e_{\max}},
\end{equation*}
where
${J^e = \det(\d \V x^e / \d\V\xi)}$.
While the aspect ratio is not affected by $p$-refinement, we observe an 
increase in level 5 and 6 due to the subdivision of deformed elements.
On the other hand, the element quality improves with each level, 
irrespective of the type of refinement.

\begin{figure}
  \subcaptionbox{Runtime of global and local cases
    \label{fig:num-exp:robustness:blade:runtimes}}
    {\includegraphics[scale=0.52]{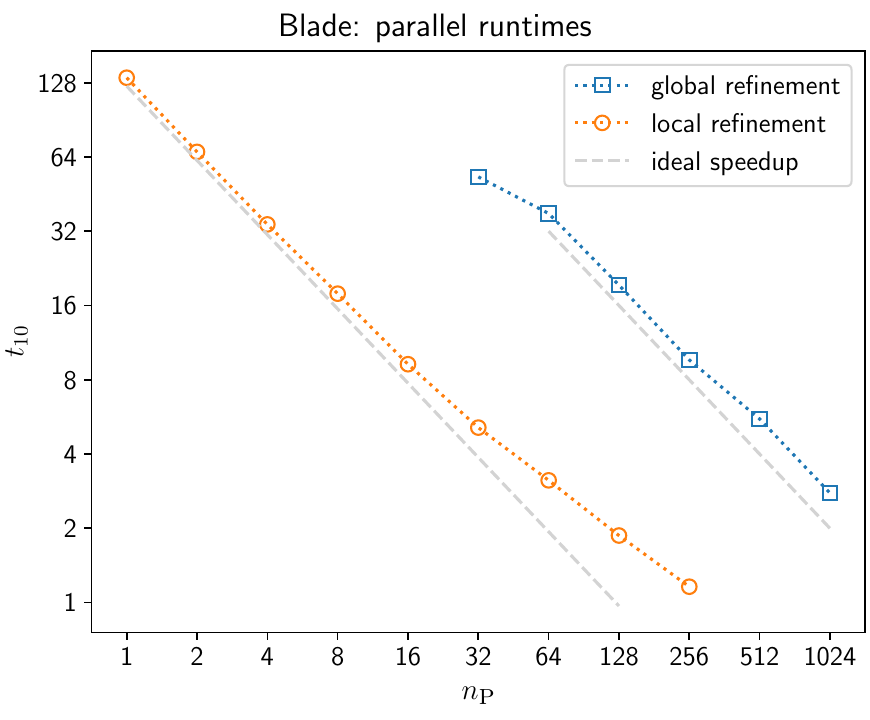}}
  \hfill
  \subcaptionbox{Imbalance of individual mesh levels in the local case
    \label{fig:num-exp:robustness:blade:imbalance}}
    {\includegraphics[scale=0.52]{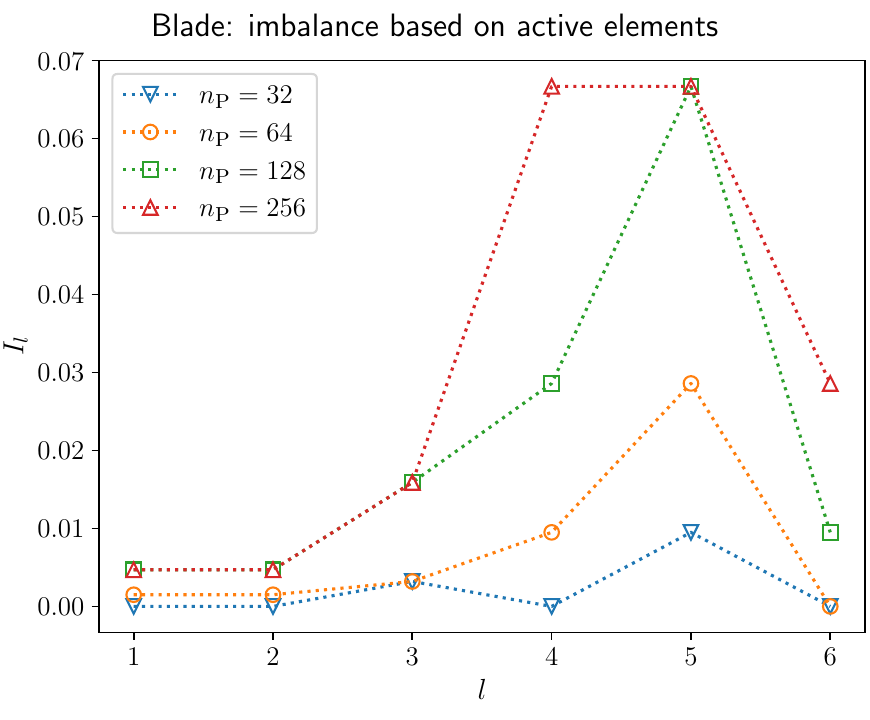}}
  \caption{Parallel performance of the FAS-FMG-SCG(6,2,4) method with
    globally and locally refined meshes.
    \label{fig:num-exp:robustness:blade}}
\end{figure}

The discrete problems were solved using FAS-FMG-SCG(6,2,4).
In the global case, ${\n{p} = 32}$ to $1024$ MPI processes were used,
and ${\n{p} = 1}$ to $256$ in the local case.
Table~\ref{tab:num-exp:robustness:blade:results} compares the
main parameters and the runtimes using 32 processes.
In spite of the challenging mesh properties, the method succeeds to reduce 
the residual by 10 orders of magnitudes with just a iteration.
Compared to the global case, the local refinement saves about 
93 percent of the memory and 90 percent of the runtime.
Per degree of freedom, the global approach requires 
$1.54\,\mu\mathrm s$ 
and the local approach
$2.12\,\mu\mathrm s$.
This is well within the range of the results achieved with slightly anisotropic
Cartesian meshes (compare Tables~\ref{tab:num-exp:robustness:cart-uni:fastest}
and \ref{tab:num-exp:robustness:cart-var:fastest}).
Figure~\ref{fig:num-exp:robustness:blade} examines the parallel performance
of the FAS-FMG-SCG(6,2,4) method with global and local refinement.
In the global case, the method exhibits an excellent strong scalability
over the full investigated range.
With local refinement, the parallel efficiency slightly degrades when using
32 or more processes.
However, this is to be expected, since as the number of processes increases, 
the partitions contain fewer active elements, while the relative proportion 
of ghost elements rises.
For example, with 256 processes, the partitions of level 4 contain just 
13 to 14 active elements.
Nevertheless, all levels remain well balanced
(Figure~\ref{fig:num-exp:robustness:blade:imbalance}).
Due to the SFC-based partitioning defined on the parent level, the number 
of active elements in the child partitions differs by no more than 8.
However, it should be noted that the workload caused by frozen elements 
was not taken into account.
Since each level is distributed across all processes, the partitions 
at the lower levels are relatively small. 
We therefore expect that the parallel efficiency can be improved by 
optimizing the partitioning strategy.

\definecolor{color3}{rgb}{1.000, 1.000, 0.498}
\definecolor{color4}{rgb}{0.667, 1.000, 0.498}
\definecolor{color5}{rgb}{0.667, 1.000, 1.000}
\definecolor{color6}{rgb}{1.000, 0.667, 1.000}

\begin{figure}
  \subcaptionbox{initial single-level mesh
    \label{fig:num-exp:robustness:blade:mesh:coarse}}
    {\includegraphics[width=0.5\textwidth]{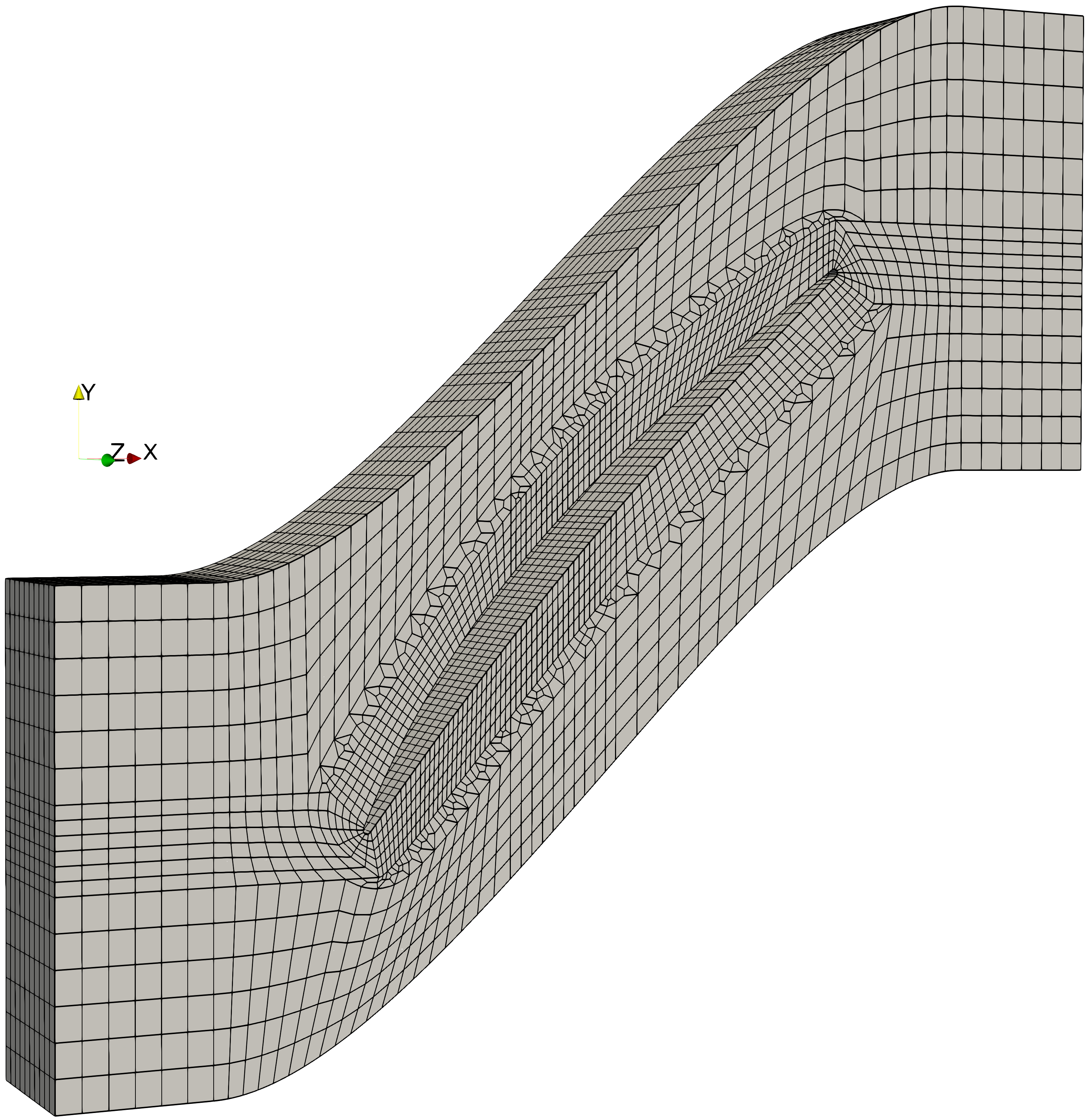}}
  \hfill
  \subcaptionbox{
    locally refined mesh near the blade nose
    \label{fig:num-exp:robustness:blade:mesh:refined}}
    {\includegraphics[width=0.45\textwidth]{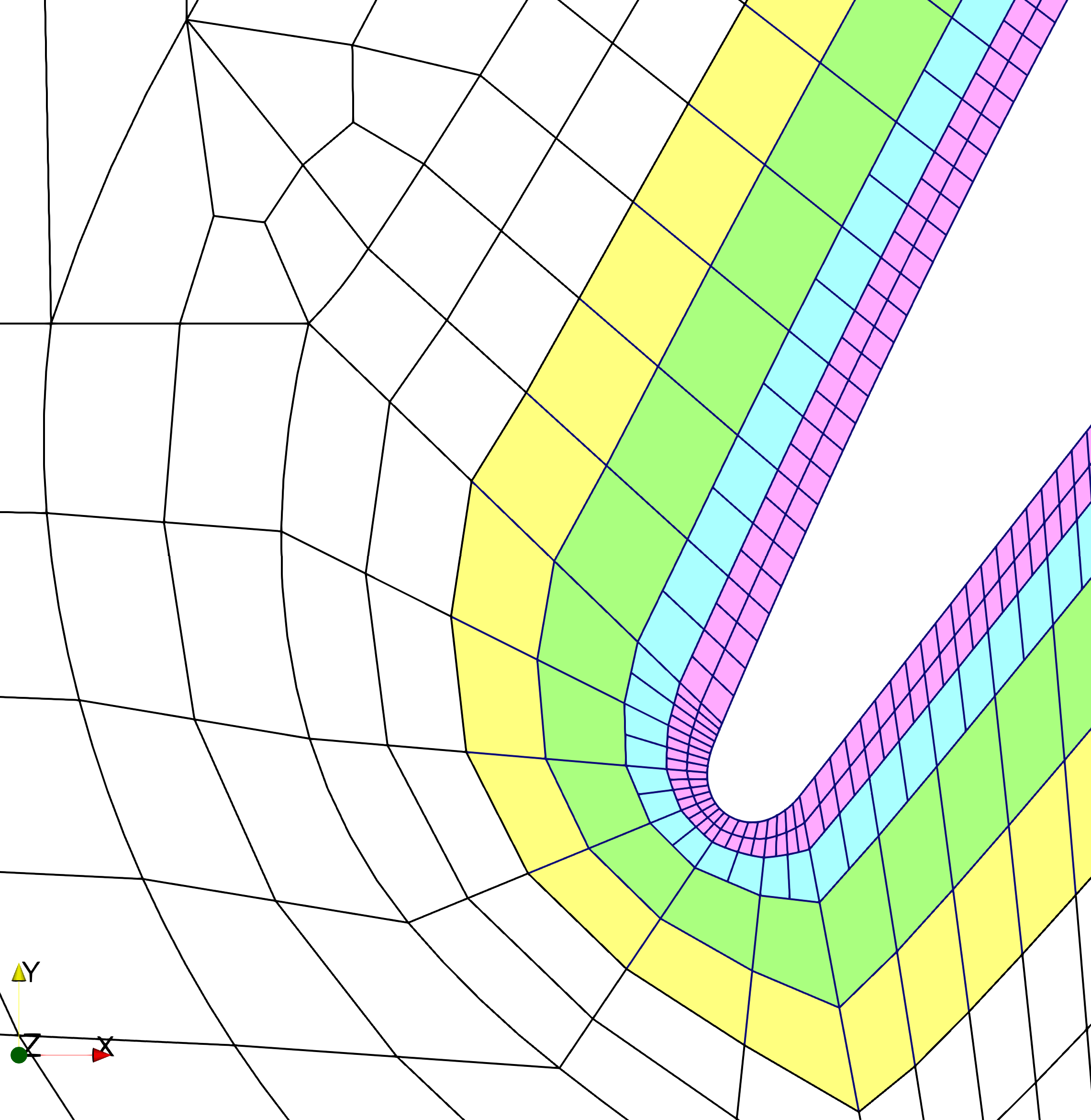}
    \vspace*{1\baselineskip}}

  \caption{Mesh around a blade in a linear compressor cascade, 
    left: coarse mesh consisting of 20640 elements;
    right: detail of 6-level composite mesh with $p/h$-refinement.
    Colors indicate the mesh levels:
    \raisebox{0.6ex}{\fcolorbox{black}{white}{\hspace*{0.5em}}}
    ${l = 2}$, ${p = 2}$;
    \raisebox{0.6ex}{\fcolorbox{black}{color3}{\hspace*{0.5em}}}
    ${l = 3}$, ${p = 4}$;
    \raisebox{0.6ex}{\fcolorbox{black}{color4}{\hspace*{0.5em}}}
    ${l = 4}$, ${p = 8}$;
    \raisebox{0.6ex}{\fcolorbox{black}{color5}{\hspace*{0.5em}}}
    ${l = 5}$, ${p = 8}$;
    \raisebox{0.6ex}{\fcolorbox{black}{color6}{\hspace*{0.5em}}}
    ${l = 6}$, ${p = 8}$.
    The subelement meshes are omitted for clarity.
    \label{fig:num-exp:robustness:blade:mesh}}
\end{figure}

\begin{table}
\footnotesize
\begin{NiceTabular}{ccrrcccc}
\toprule
 $l$ 
 & $p_l$ 
 & $\n[,l]{e}^\mathsc{a}$ 
 & $\n[,l]{e}^\mathsc{l}$ 
 & $\Delta x^e_{\min}$ 
 & $\Delta x^e_{\max}$ 
 & $\mathit{AR}_{l,\max}$
 & $Q_{l,\min}$
 \\
\midrule

  1  &  1  &  20640  &      0   &  8.751E-03  &  3.370E-02  &  5.827  &  0.151  \\
  2  &  2  &  20640  &  15600   &  8.751E-03  &  3.370E-02  &  5.827  &  0.144  \\
  3  &  4  &   5040  &   1680   &  9.116E-03  &  2.014E-02  &  5.827  &  0.364  \\
  4  &  8  &   3360  &   1680   &  9.116E-03  &  2.014E-02  &  5.827  &  0.364  \\
  5  &  8  &  13440  &   6720   &  4.209E-03  &  1.018E-02  &  7.395  &  0.522  \\
  6  &  8  &  53760  &  53760   &  1.998E-03  &  5.119E-03  &  8.208  &  0.685  \\ 

\bottomrule
\end{NiceTabular}

\caption{Characteristics of the locally refined multilevel mesh for the
   compressor cascade.}
\label{tab:num-exp:robustness:blade:mesh-lr}
\end{table}

\begin{table}
\footnotesize
\begin{NiceTabular}{l%
                    @{\quad}>{\raggedleft\arraybackslash}p{5em}%
                    @{\quad}>{\raggedleft\arraybackslash}p{5em}}
\toprule
     & global& local \\
\midrule 

$\n{dof}^{\mathrm{tot}}$  &  1102\,M   &  79\,M  \\
$\n{dof}^{\mathsc{l}}$    &   963\,M   &  46\,M  \\
$n_{10}$                  &        1   &      1  \\
$t_{10}(32)$              &    53.14   &   5.22  \\
$\tau_{10}(32)$           &     1.54   &   2.12  \\

\bottomrule
\end{NiceTabular}

\caption{Parameters and performance metrics of FAS-FMG-SCG(6,2,4) 
  for blade test case with global and local refinement, respectively.}
\label{tab:num-exp:robustness:blade:results}
\end{table}

\subsection{Parallel adaptive multigrid}
\label{sec:num-exp:adaptive}

\begin{figure}
  \includegraphics[scale=0.1]{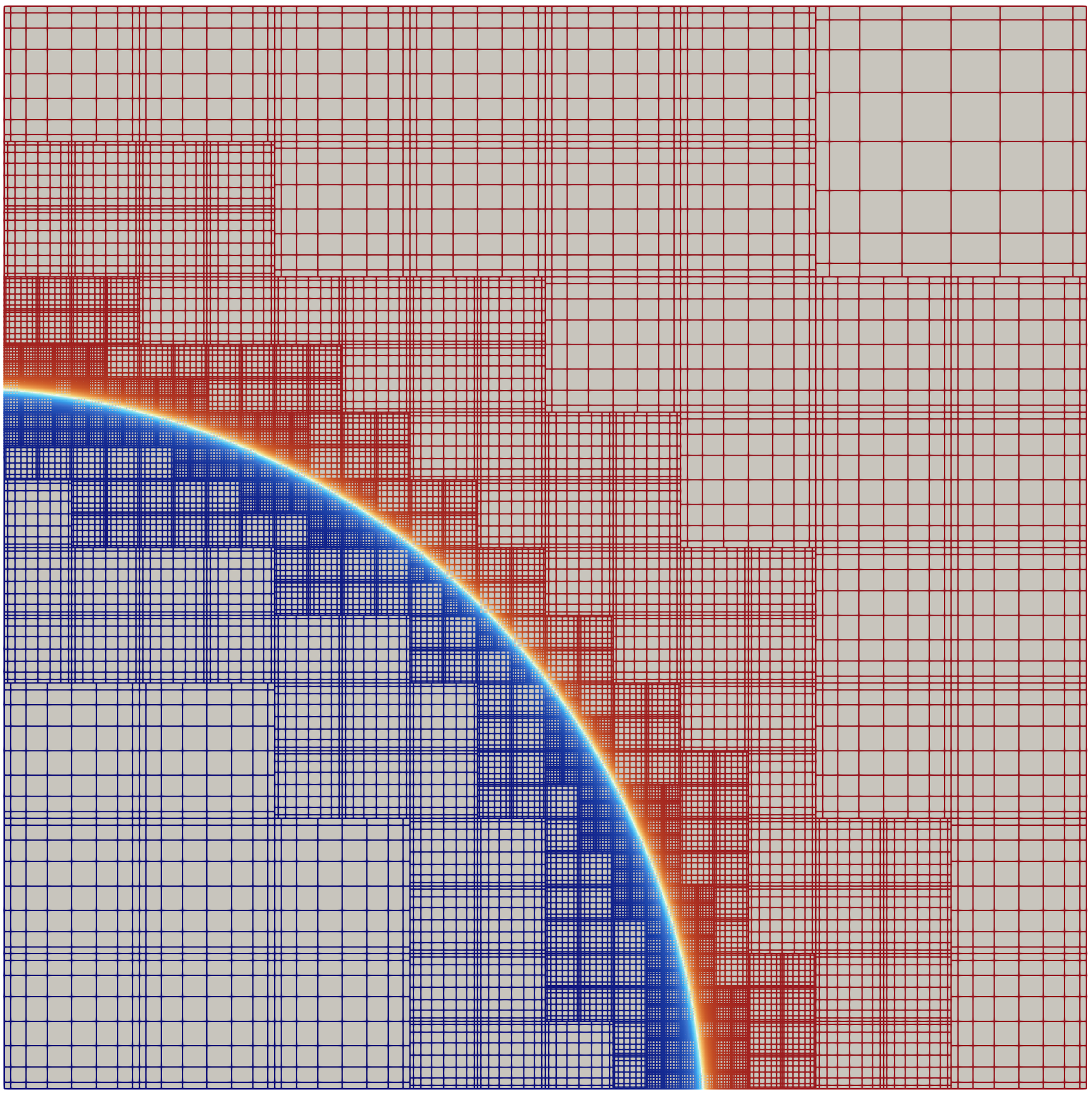}
  \caption{Side view of the composite adapted mesh with 10 levels and $p=8$,
    coloring based on the numerical solution.
    \label{fig:num-exp:adaptive:mesh}}
\end{figure}

In the previous test cases, the focus was on robustness and the fast solution of
discrete equations.
To investigate the potential savings from mesh adaptation, we use the  FAS-FMG
method to solve the benchmark problem proposed in \cite{MM_Cerveny2019a}.
The problem is equivalent to \eqref{eq:helmholtz:pde} on the unit cube 
${\Omega = (0,1)^3}$
with ${\kappa = 0}$
and Dirichlet boundary conditions.
The right-hand side is chosen so that the exact solution is
\begin{equation*}
  u(\V x) = \arctan\left(\alpha (\abs{\V x - \V x_{\mathrm c}} - r) \right)
  \,.
\end{equation*}
The solution is characterized by a sharp spherical wave front of radius $r$
centered at ${\V x_{\mathrm c}}$.
Following \cite{MM_Cerveny2019a} we choose
${\alpha = 200}$,
${\V x_{\mathrm c} = -[0.05,0.05,0.05]}$ and
${r = 0.7}$.
The initial mesh is constructed starting from the root level with 
$4^3$ cubic elements with ${p_1 = 1}$ by applying a sequence of global
$p$-refinements, until the targeted polynomial degree is reached on 
level $L_0$.
Then, an adaptation loop is executed, comprising the following steps:
\begin{enumerate}
\item
  Solve the problem on the current mesh using two cycles of FAS-FMG-SCG(2,2,2)
\item
  For each level ${l \ge L_0}$, compute the energy norm of every element
  $\Omega_l^e$
  \begin{equation*}
    \varepsilon_{l,\mathsc e}^{e}
      = \norm{\varepsilon_l}_{\mathsc e}^{e}
      = \NM{\varepsilon}_{l}^{e} 
        \NM A_{l}^{e} \NM{\varepsilon}_{l}
  \end{equation*}
\item
  Record the total error 
  $\varepsilon_{\mathsc e}$
  and the number of leaf elements in each level
\item
  Refine leaf elements in ${l < L_{\max}}$ for which
  ${
    \varepsilon_{l,\mathsc e}^{e} > 
    0.7 \max_{l,e}(\varepsilon_{l,\mathsc e}^{e})
  }$
\end{enumerate}
To calculate the element error norm in step 2, 
the discrete operator is first applied to the complete error vector, 
then restricted to the element, 
and finally weighted using the element error coefficients.
This discrete energy norm accounts for the jumps across the element boundary,
which are absent in the continuous FEM applied in \cite{MM_Cerveny2019a}.
While this norm is asymptotically equivalent to its continuous counterpart, 
it tends to underestimate the error if the resolution is too low.
To compensate for this, the elements cut by the wavefront are refined 
regardless of the error.
Please note that the number of levels is restricted to ${L = L_{\max}}$.
Elements in this level are not further refined.
Figure~\ref{fig:num-exp:adaptive:mesh} shows an adapted mesh of degree 8
with ten levels.
Note that the refinement extends beyond the immediate vicinity of the front, 
since distant elements contribute progressively to the error as the accuracy
increases.

\begin{figure}
  \subcaptionbox{Convergence with global and adaptive refinement
    \label{fig:num-exp:adaptive:error}}
    {\includegraphics[scale=0.52]{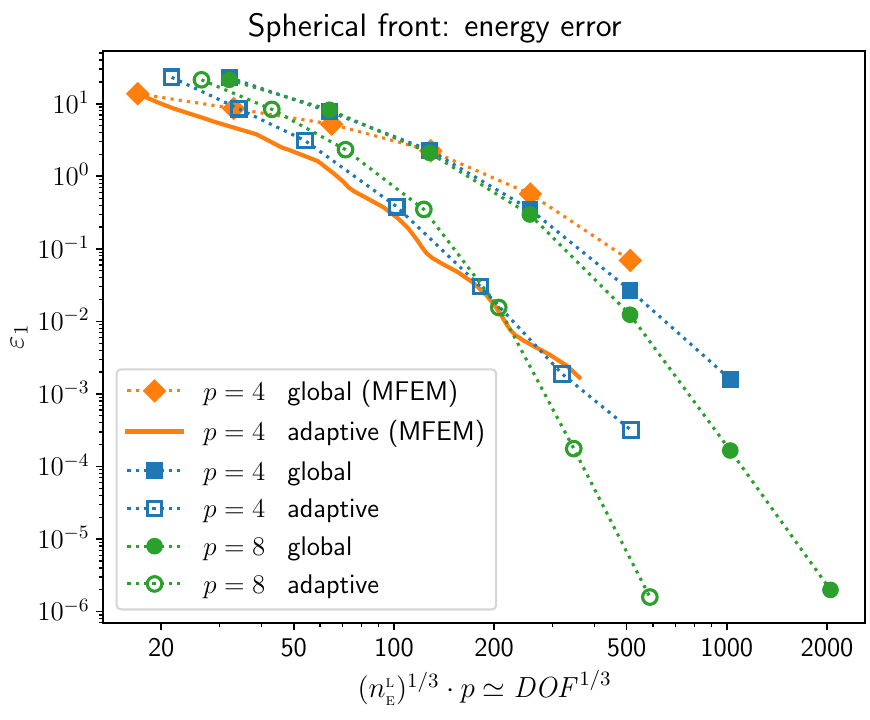}}
  \hfill
  \subcaptionbox{Normalized runtimes of substeps
    \label{fig:num-exp:adaptive:runtimes}}
    {\includegraphics[scale=0.52]{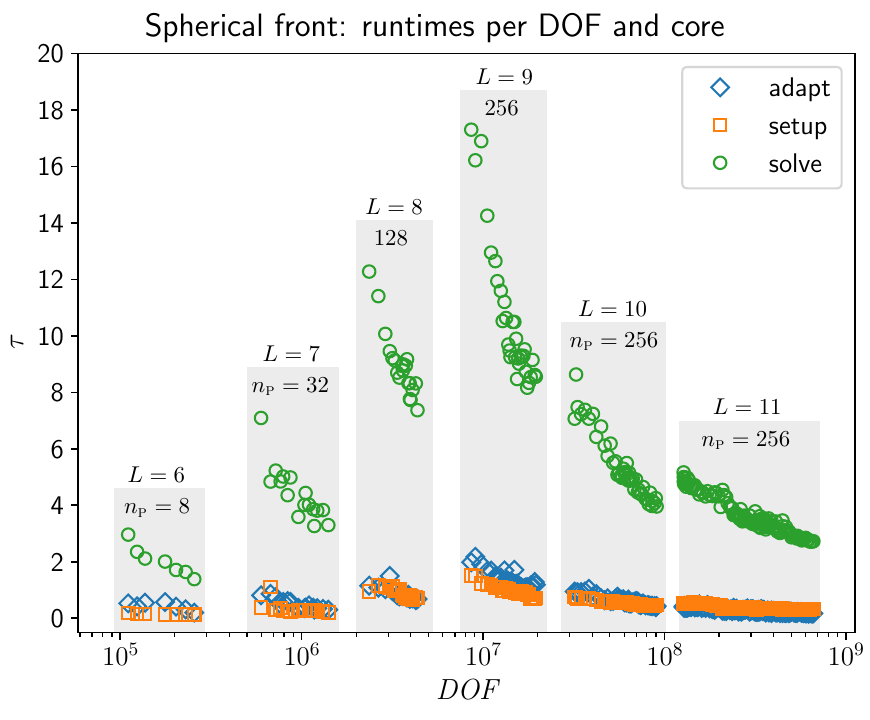}}
  \caption{Spherical front test case. Left: convergence of adaptive FAS-FMG 
    compared to global refinement and MFEM results from \cite{MM_Cerveny2019a}
    (courtesy Jakub \v{C}erven\'{y}). 
    Right: normalized runtimes of the adaptation, operator-setup and solution 
    steps using a different number of mesh levels $L$ and MPI processes $\n{p}$.
    Every symbol represents a single iteration of the adaptation loop.
    \label{fig:num-exp:adaptive:results}}
\end{figure}

Figure~\ref{fig:num-exp:adaptive:error} compares the convergence of
the adaptive FAS-FMG method for ${p=4}$ and ${p=8}$ with corresponding
results using global refinement.
For ${p = 4}$, our results agree well with those obtained using the 
continuous finite element method implemented in MFEM 
\cite{MM_Cerveny2019a}.
Using a higher approximation order (${p = 8}$) is initially slightly 
more expensive, but becomes substantially more efficient for errors
below $10^{-2}$.
In both cases, between 93\% and 98\% DOFs are saved compared to global 
refinement (increasing with decreasing error).
Figure~\ref{fig:num-exp:adaptive:runtimes} compares the runtimes
of the mesh adaptation and solution steps using a different number of mesh 
levels and MPI processes.
Each shaded box indicates a separate adaptation loop.
Please note that every loop starts with the initial multilevel mesh, 
which contains only a few elements in each partition.
This leads to scaled solver times that are significantly greater than in 
the previous tests using a static multilevel mesh.
However, as the loop progresses, the number of element grows and the solver
time drops into the range between 1 and 4 microseconds per DOF and core. 
It is also worth noting that mesh adaptation and operator setup both require only
about 10\% of the solver time.
Altogether these results demonstrate the excellent efficiency and parallel 
scalability of our method within the investigated parameter range.

\subsection{Comparison to the state-of-the art}
\label{sec:num-exp:state-of-the-art}

The multigrid methods developed in this work differ from common approaches 
in two respects:
\begin{enumerate}
\item
  the application of Krylov-accelerated weighted Schwarz smoothers, and
\item
  the use of the FAS scheme on locally refined grids.
\end{enumerate}
Recently, polynomial smoothers such as the Chebyshev method gained 
considerable attention, e.g. 
\cite{MG_Fehn2020a,MG_Munch2023b,MG_Phillips2025a}.
\citet{MG_Fehn2020a} conducted extensive performance tests with Chebyshev-Jacobi
smoothers using a single core of an Intel Xeon 8174. 
This CPU has a $1.35$ times higher peak performance per core than 
the Intel Xeon 8470, but lower memory bandwidth.
Nevertheless, the results should be roughly comparable with the present study.
For Cartesian meshes, the authors reported runtimes ranging from
${\tau_{10} = 2.5}$ for ${p = 3}$ to
${\tau_{10} = 7}$ for ${p = 15}$.
In contrast, the fastest methods in the current study achieved 
${\tau_{10} = 1.8}$ for ${p = 4}$ and
${\tau_{10} = 0.3}$ for ${p = 8}$ and $16$.
While the performance is comparable at lower degrees, our results indicate a
speedup by a factor of 20 or more for ${p \ge 8}$.
Given the similar hardware characteristics, this improvement is most likely 
attributed to the higher algorithmic efficiency of the present approach.
In a follow-up study, \citet{MG_Munch2023b} developed Schwarz-Chebyshev 
smoothers serving as a multigrid preconditioner for flexible CG or GMRES 
methods.
For Cartesian grids with an aspect ratio of ${\mathit{A} = 50}$, 
the computational cost increased by a factor of 10 compared to 
the isotropic case.
Similar observations were made in \cite{MG_Phillips2025a}
In contrast, the Krylov-accelerated Schwarz smoothers developed in the present study 
yield efficient multigrid methods that are practically robust against the aspect ratio.
For instance, consider FAS-FMG-SCG(6,2,4), as shown in 
Figure~\ref{fig:num-exp:cart-var}.

Unlike other recently proposed high-order adaptive multigrid methods,
our method does not operate on a global adaptive mesh, but on a locally
refined multilevel mesh.
Each mesh level extends only over the corresponding refinement zone.
This approach avoids the introduction of hanging nodes and related
modifications of the discrete operators.
Starting from the root mesh, all mesh levels result naturally from the refinement 
process and need not to be constructed by coarsening the global top level mesh.
For parallelization, the meshes are partitioned level by level.
By using a recursively refined space-filling curve, the partitioning requires 
only a fraction of the execution time of a single V-cycle and limits the 
imbalance to no more than 8 elements on each level.
This strategy will most likely produce smaller imbalances across the levels 
than the global coarsening and local smoothing methods proposed in
\cite{MG_Clevenger2020a,MG_Munch2023a}.
However, these and further performance aspects are beyond the scope of the 
present study and will be addressed in future work.



\section{Prospects for application to flow problems}
\label{sec:outlook}


%
Many time-integration methods for the Navier-Stokes equations result in elliptic subproblems that must be solved at each time step. 
For example, the second-order velocity correction scheme for incompressible flow 
reads \cite{TI_Guermond2003b,TI_Guesmi2023a}
\begin{subequations}
\label{eq:ins:v-correction}
\begin{gather}
  \label{eq:ins:extrapolation}
  \frac{\gamma_0 \V v' - \alpha_0 \V v^{n} - \alpha_1 \V v^{n-1}}{\Delta t}
  = -\sum_{i=0}^{1} 
       \beta_{i} \bigl( \nabla\cdot(\V v \V v) 
                      + \nu\nabla\times\nabla\times\V v \bigr)^{n-i}
\\
  \label{eq:ins:projection}
  \frac{\gamma_0 (\V v" - \V v')}{\Delta t} = - \nabla p^{n+1}
  , \quad 
  \nabla \cdot \V v" = 0
\\
  \label{eq:ins:diffusion}
  \frac{\gamma_0}{\Delta t} \V v^{n+1}
  - \nabla\cdot\nu\bigl(\nabla\V v + \transpose{(\nabla\V v)}\bigr)^{n+1}
  = \frac{\gamma_0}{\Delta t} \V v"
  + \sum_{i=0}^{1} 
     \beta_{i} \nu\nabla\times\nabla\times\V v^{n-i},
\end{gather}
\end{subequations}
where
$\V v$ denotes the velocity,
$p$ the pressure and
$\nu$ the kinematic viscosity.
Note that the equations are written in dimensionless form, so that the density does not appear.
Further, 
$n$ indicates the time step,
$\Delta t$ is the step size and
${\alpha_0 = 2}$,
${\alpha_1 = -1/2}$,
${\beta_0  = 2}$,
${\beta_1  = -1}$,
${\gamma_0 = 3/2}$
are the coefficients of the underlying semi-implicit backward difference formula.
The scheme begins with an explicit extrapolation in Equation \eqref{eq:ins:extrapolation}.
Following this, the projection \eqref{eq:ins:projection} eliminates the divergence 
from the velocity.
As an intermediate step, the divergence of the first sub-equation is taken and combined with the second to obtain the pressure equation
\begin{equation}
  \label{eq:ins:pressure}
  -\nabla^2 p^{n+1} = -\frac{\gamma_0}{\Delta t} \nabla \cdot \V v'
  \,.
\end{equation}
This equation yields the pressure $p^{n+1}$, which is employed to project
the extrapolated velocity $\V v'$ to the solenoidal field $\V v"$.
Solving the viscous diffusion problem \eqref{eq:ins:diffusion} completes
the time step.
Equations \eqref{eq:ins:v-correction} and \eqref{eq:ins:pressure} are discretized using 
an FAS-multilevel DG method with local Lax-Friedrichs fluxes for convection.
Following the extrapolation step, the discrete pressure equation is solved by 
means of the FAS multigrid method developed in Section~\ref{sec:mg-method}.
Please note that no startup procedure is required, as reasonable initial values 
are available from the previous time step.
Finally, the diffusion system is solved using FAS-MG with a slightly modified 
SCG smoother, which applies the Schwarz preconditioner to each component 
individually.  

\begin{figure}
  \includegraphics[width=\textwidth]{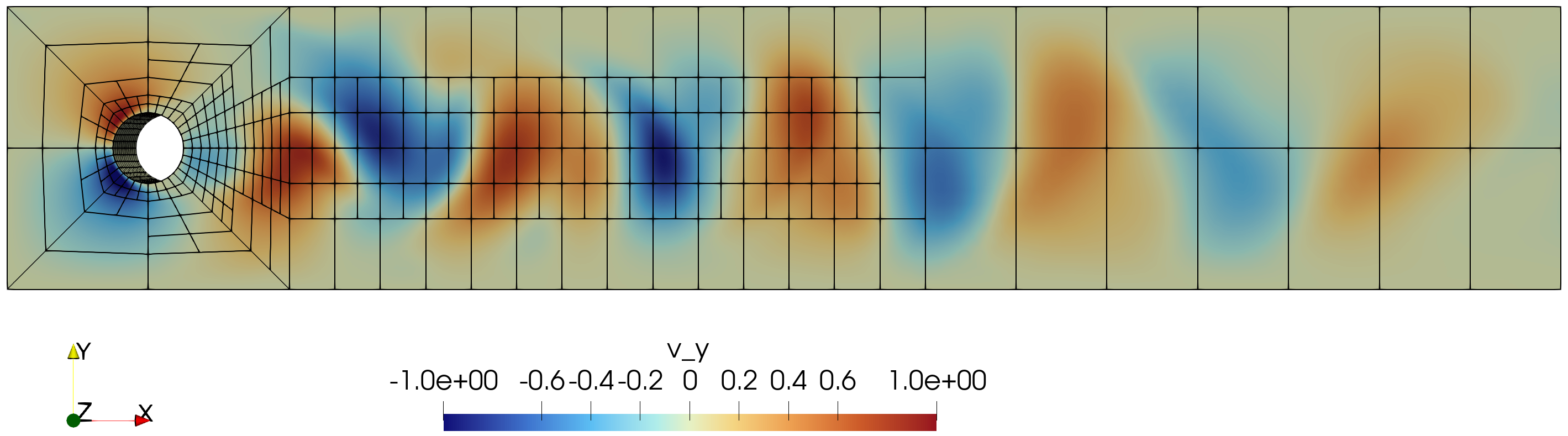}
  \caption{Preliminary results for the cylinder flow benchmark 
    \cite{BM_John2004a}: composite element mesh and $y$-component of velocity at
    ${t=5}$.
    \label{fig:outlook:cylinder-mesh-vy}}
\end{figure}

Figure~\ref{fig:outlook:cylinder-mesh-vy}
shows preliminary results of the method for the unsteady flow past a circular cylinder 
at Reynolds number $100$, proposed as a benchmark problem in \cite{BM_John2004a}.
For this study, the problem was discretized using five mesh levels with polynomial 
degrees ${\NM p = \{ 2,3,7,7,7 \}}$.
Starting with 156 elements on the root level, the mesh was locally refined near 
the cylinder and in the wake on levels 4 and 5.
The final composite mesh comprises 3096 elements of degree $7$, 
yielding $6.34$ million degrees of freedom for velocity and pressure.
Thorough investigations, application to turbulent flows and the extension to 
higher-order time-integration methods are the subject of ongoing work.

\section{Conclusion}
\label{sec:conclusion}

In this work, we develop a multigrid method for solving elliptic problems that are discretized using the discontinuous Galerkin spectral element method on hexahedral meshes.
The method extends the full-approximation storage algorithm of \citet{MG_Brandt1977a} by introducing frozen elements that allow for seamless coupling across refinement zones without the need for hanging nodes.
Unlike common approaches, the solution is stored on multiple mesh levels instead of being defined on a locally refined mesh that covers the entire computational domain.
$\mathit{hp}$-adaptivity is achieved by dynamically adjusting the resulting multilevel mesh.

The core component of the multigrid solver is the smoother based on an overlapping Schwarz method, which has been generalized from \cite{MG_Stiller2017a} to handle unstructured curvilinear meshes. This smoother preserves the tensor product structure, allowing the application of the fast diagonalization method presented in \cite{MG_Lynch1964a}.
To enhance robustness, the Schwarz method is accelerated using the inexact conjugate gradient method of \citet{KR_Golub1999a}.
By combining this approach with the full multigrid startup strategy, we obtain an FAS-FMG method tailored for DG-SEM on locally refined meshes.
Numerical experiments demonstrate the exceptional efficiency of this method and its robustness against high aspect ratios, element deformation, and irregular mesh topology.
Although a direct comparison is difficult, the results suggest that our method surpasses recently proposed MG methods based on Chebyshev smoothers and other polynomial smoothers \cite{MG_Fehn2020a,MG_Munch2023a,MG_Phillips2025a}.

In this study, we focused on the algorithmic aspects.
Our ongoing work is devoted to investigating and improving parallel efficiency and scalability of the method.
Future research will expand the method to more intricate applications, such as flow problems, and combine it with error estimation methods for guiding mesh adaptation


\bibliographystyle{abbrvnat}
\bibliography{BM,FE,FV,HPC,KR,MG,MM,SE,TB,TI}




\end{document}